\documentclass{amsart} 
\usepackage{graphicx}
\usepackage{amssymb, amsmath, sidecap}
\usepackage{amsfonts}
\usepackage{amssymb}
\usepackage{float}
\usepackage{extarrows}
\usepackage{mathrsfs}
\usepackage{booktabs}
\usepackage{verbatim}
\usepackage{hyperref}
\usepackage{esint}
\usepackage[usenames,dvipsnames]{xcolor}

\renewcommand{\baselinestretch}{1}

\newcommand{\dist}{\text{\rm dist}}

\def\bt{\begin{thm}}
\def\et{\end{thm}}
\def\bl{\begin{lem}}
\def\el{\end{lem}}
\def\bd{\begin{defi}}
\def\ed{\end{defi}}
\def\bc{\begin{cor}}
\def\ec{\end{cor}}
\def\bp{\begin{proof}}
\def\ep{\end{proof}}
\def\br{\begin{rem}}
\def\er{\end{rem}}

\def\Forall{\text{ } \forall \:}
\def\d{\, \mathrm{d}}

\def\be{\begin{equation}}
\def\ee{\end{equation}}
\def\bes{\begin{equation*}}
\def\ees{\end{equation*}}
\def\bea{\begin{equation} \begin{aligned}}
\def\eea{\end{aligned} \end{equation}}
\def\beas{\begin{equation*} \begin{aligned}}
\def\eeas{\end{aligned} \end{equation*}}
\def\ba{\begin{align}}
\def\ea{\end{align}}
\def\bas{\begin{align*}}
\def\eas{\end{align*}}

\newtheorem{thm}{Theorem}[section]
\newtheorem{lem}{Lemma}[section]
\newtheorem{defi}{Definition}[section]

\newtheorem{prop}[thm]{Proposition}

\newtheorem{rem}{Remark}[section]
\newtheorem{cor}{Corollary}[section]

\numberwithin{equation}{section}
\numberwithin{figure}{section}

\begin{document}
\title{ $L^{\Phi(A(x,\cdot))}$ estimate for the gradient in $W^{1,A(x,\cdot)}$ \footnotemark[1]\footnotemark[2]}

\author[Liu]{Duchao Liu}
\address[Duchao Liu]{School of Mathematics and Statistics, Lanzhou
University, Lanzhou 730000, P. R. China} \email{liuduchao@gmail.com,
Tel.: +8613893289235, fax: +8609318912481}

\author[Wang]{Beibei Wang}
\address[Beibei Wang]{School of Mathematics and Statistics, Lanzhou
University, Lanzhou 730000, P. R. China} \email{wangbb15@lzu.edu.cn}

\author[Yang]{Sibei Yang}
\address[Sibei Yang]{School of Mathematics and Statistics, Lanzhou
University, Lanzhou 730000, P. R. China} \email{yangsb@lzu.edu.cn}

\footnotetext[2]{
Research supported by the National Natural Science Foundation of
China (NSFC 11501268 and NSFC 11471147).
}

\begin{abstract}
Under appropriate assumptions on the $N(\Omega)$-fucntion, the $L^{\Phi(A(x,\cdot))}$ estimate for the gradient of the minimizers of a class of energy functional in Musielak-Orlicz-Sobolev space $W^{1,A(x,\cdot)}$ is presented by using Calder\'{o}n-Zygmund decomposition.
\end{abstract}
 \maketitle


\section{Introduction}

Vast mathematical literature describes various aspects of partial differential equations related to the elliptic type operators including variable exponent, weighted, convex and double phase cases. Musielak-Orlicz-Sobolev spaces give an abstract framework of functional analysis to cover all of the above mentioned cases.

Some basic references of special Musielak-Orlicz spaces should be noticed. Before the key role in the modular spaces of functional analysis was played by the comprehensive book by Musielak \cite{Musielak}, Nakano \cite{Nakano50} provided the first framework approach to non-homogeneous setting with general growth conditions, which was followed by Skaff \cite{Skaff69I, Skaff69II} and Hudzik \cite{Hudzik76I, Hudzik76II}. In much earlier time variable exponent spaces were introduced by the pioneering work by Orlicz \cite{Orlicz31}.

\subsection{Linear and nonlinear Calder\'on-Zygmund theory}

In this paper we mainly concern with the gradient integrability properties of the minimizers for some energy functional in Musielak-Orlicz-Sobolev space $W^{1,A(x,\cdot)}$ under some reasonable assumptions on the $N(\Omega)$-function $A$. Our results can be basically considered as the generalization of the classical Calder\'on-Zygmund theory for linear case started by the work of Calder\'on and Zygmund \cite{CZ52, CZ56}. From the viewpoint of contemporary mathematical research, these two foundation papers \cite{CZ52, CZ56} by Calder\'on and Zygmund yield a new research field in mathematics - the theory of harmonic analysis. Since the linear Calder\'on-Zygmund theory heavily depends on the existence of the fundamental solution, the nonlinear case of the theory does not exist before De Giorgi's famous iteration method firstly appeared in \cite{DeGiorgi57}. 

The nonlinear case of the Calder\'on-Zygmund theory was started  in the fundamental work of Iwaniec \cite{Iwaniec82, Iwaniec83}. Then important contributions by DiBenedetto and Manfredi \cite{DiBenedetto93} and Caffarelli and Peral \cite{Caffarelli98} followed the time line of the investigation. After those papers, a large number of literature on the Nolinear Calder\'on-Zygmund theory appears in recent years, see for instance \cite{Kinnunen99, Gutierrez01, Byun09, Byun12, Mingione07, Mingione10, Mingione15}.




\subsection{Recent regularity results in Musielak-Orlicz-Sobolev spaces}

A highly important part of the mathematical literature in general Musielak-Orlicz-Sobolev spaces gives structural conditions on regularity analysis of the space in recent years.  

In a recent work \cite{Ahmida18}, Ahmida and collaborators prove the density of smooth functions in the modular topology in Musielak-Orlicz-Sobolev spaces, which extends the results of Gossez \cite{Gossez82} obtained in the Orlicz-Sobolev settings. The authors impose new systematic regularity assumption on the modular function. And this allows to study the problem of density unifying and improving the known results in Orlicz-Sobolev spaces, as well as variable exponent Sobolev spaces. 

In the paper \cite{liu_yao_18}, under some reasonable assumptions on the $N(\Omega)$-fucntion, the De Giorgi process is presented by the authors in the framework of Musielak-Orlicz-Sobolev spaces. And as the applications, the local bounded property of the minimizer for a class of the energy functional in Musielak-Orlicz-Sobolev spaces is proved. Under similar assumptions as in \cite{liu_yao_18}, the authors in \cite{Liu18} prove the
H\"{o}lder continuity of the minimizers for a class of the energy
functionals in Musielak-Orlicz-Sobolev spaces.

\subsection{Motivation}

The priori estimate of $L^{\Phi(A(x,\cdot))}$ for the gradient in $W^{1,A(x,\cdot)}$ relates to many aspect of regularity results in the Musilak-Orlicz-Sobolev spaces, for example $C^{1,\alpha}$ estimate. Similar priori estimate results in the classical Sobolev spaces can be found in \cite{Bojarski55, Gehring73, Gianquinta12, Giaquinta79}. In the paper \cite{Fan96}, Fan presents the $L^{q(\cdot)}$ and $L^\infty$ estimate for the gradient of the minimizers of some kind of functionals in $W^{1,p(\cdot)}$. In this paper we generalized the results of Fan. We emphasize that the assumptions to derive the regularity results in this paper cover the variable exponent case, but can not cover the ones in double phase case in \cite{Mingione15}. In order to get the estimate of $L^{\Phi(A(x,\cdot))}$ for the gradient, it is still interesting to find out the uniform assumptions including both variable exponent case and double phase case in Musielak-Orlicz-Sobolev spaces in order to derive these regularity results.

To find a more reasonable method to describe increasing conditions for the modular functions is always hard to overcome. Moreover, we wish the method we have developed are convenient to make the calculus for $N(\Omega)$ or $N$ functions. In this paper, based on the tools developed in \cite{liu_yao_18, Liu18} (Lemma \ref{compare}), we find some additional more accurate ways to describe increasing conditions on the $N(\Omega)$ functions in order to make calculus of ``functions'' accurately, see $\Delta_{\mathbb{R}^+}$ and $\nabla_{\mathbb{R}^+}$ conditions in Definition \ref{DelT} and their basic properties in Lemma \ref{bac_delta1} and Lemma \ref{bac_delta2}. For these basic assumptions, our fundamental viewpoint is that the more accurate estimate needs more detailed assumptions. 

In most of existing research papers, the mathematicians are used to describing increasing condition by power functions. But it is natural to do this by the modular functions in Musielak-Orlicz spaces instead of the power functions. Furthermore, it is natural to describe modular functions by the convexity of the function in Musielak-Orlicz spaces. To see this, for example, it is well known that the H\"{o}lder inequality holds for convex modular functions (see Proposition \ref{int_A} (5)). We will use the uniform convexity to describe the increasing condition on modular functions in this paper.

To overcome the linearity dependence in the linear Calder\'on-Zygmund theory, we establish the Caccioppoli type inequality by De Giorgi's iteration, see Lemma \ref{ke_inl}. And by a new type of Gehring type inequality developed in Lemma \ref{ke_ineq} and an abstract theorem developed by Calder\'on-Zygmund decomposition in Theorem \ref{kkeeyyLemma}, we rebuild the nonlinear $L^{\Phi(A(x,\cdot))}$ estimate for the gradient in $W^{1,A(x,\cdot)}$.

\vspace{0.4cm}
The main results of this paper are Theorem \ref{kkeeyyLemma} in Section \ref{Sec3} and Theorem \ref{main_ltheorem} in Section \ref{Sec4}. These results extend part of the results in \cite{ Gehring73, Giaquinta79, Fan96}. We also generalize the log-H\"{o}lder continuity of the variable exponent case in Musielak-Orlicz-Sobolev spaces in Section \ref{Sec4}, see Lemma \ref{LHol}.

\section{The Musielak-Orlicz-Sobolev Spaces}\label{Sec2}

In this section, we list some definitions and propositions related to
Musielak-Orlicz-Sobolev spaces.
Firstly, we give the definition of \textit{$N$-function} and
\textit{generalized $N$-function} as following.

\vspace{0.3cm}

\begin{defi}
A function $A:\mathbb{R}\rightarrow[0,+\infty)$ is called an
$N$-modular function (or $N$-function), denoted by $A\in N$, if $A$ is even and convex,
$A(0)=0, 0< A(t)\in C^0$ for $t\not=0$, and the following conditions
hold
\begin{equation*}
\lim_{t\rightarrow0+}\frac{A(t)}{t}=0\text{ and }
\lim_{t\rightarrow+\infty}\frac{A(t)}{t}=+\infty.
\end{equation*}
Let $\Omega$ be a smooth domain in $\mathbb{R}^n$. A function $A:\Omega\times\mathbb{R}\rightarrow[0,+\infty)$ is
called a generalized $N$-modular function (or generalized $N$-function), denoted by $A\in N(\Omega)$, if
for each $t\in[0,+\infty)$, the function $A(\cdot,t)$ is measurable,
and for a.e. $x\in\Omega$, we have $A(x,\cdot)\in N$.
\end{defi}

\vspace{0.3cm}

Let $A\in N(\Omega)$, the Musielak-Orlicz space $L^{A}(\Omega)$ is
defined by
\begin{equation*}
\begin{aligned}
L^{A}(\Omega)&:=\bigg\{u:\,u\text{ is a measurable real function, and
}\exists\lambda>0\\
&\quad\quad\quad\quad\quad\quad\quad\quad\quad\quad\quad\quad\text{ such that
}\int_{\Omega}A\bigg(x,\frac{|u(x)|}{\lambda}\bigg)\,\mathrm{d}x<+\infty\bigg\}
\end{aligned}
\end{equation*}
with the (Luxemburg) norm
\begin{equation*}
\|u\|_{L^{A}(\Omega)}=\|u\|_A:=\inf\bigg\{\lambda>0:\,\int_{\Omega}A\bigg(x,\frac{|u(x)|}
{\lambda}\bigg)\,\mathrm{d}x\leq1\bigg\}.
\end{equation*}

The Musielak-Sobolev space $W^{1,A}(\Omega)$ can be defined by
\begin{equation*}
W^{1,A}(\Omega):=\{u\in L^{A}(\Omega):\,|\nabla u|\in
L^{A}(\Omega)\}
\end{equation*}
with the norm
\begin{equation*}
\|u\|_{W^{1,A}(\Omega)}=\|u\|_{1,A}:=\|u\|_A+\|\nabla u\|_{A},
\end{equation*}
where $\|\nabla u\|_{A}:=\|\,|\nabla u|\,\|_{A}$.

$A$ is called locally integrable if $A(\cdot,t_0)\in
 L_{\text{loc}}^1(\Omega)$ for every $t_0>0$.

\vspace{0.3cm}

\begin{defi} We say that $a(x,t)$ is the Musielak derivative of $A(x,t)\in N(\Omega)$ at $t$ if
for $x\in\Omega$ and $t\geq0$, $a(x,t)$ is the right-hand derivative of
$A(x,\cdot)$ at $t$; and for $x\in\Omega$ and $t\leq0$,
$a(x,t):=-a(x,-t)$.
\end{defi}

\vspace{0.3cm}

Define $\widetilde A:\Omega\times\mathbb{R}\rightarrow[0,+\infty)$
by
\begin{equation*}
\widetilde A(x,s)=\sup_{t\in\mathbb{R}}\big(st-A(x,t)\big)\text{ for
}x\in\Omega\text{ and }s\in\mathbb{R}.
\end{equation*}
$\widetilde A$ is called the \textit{complementary function} to $A$ in the
sense of Young. It is well known that if $A\in N(\Omega)$, then
$\widetilde A\in N(\Omega)$ and $A$ is also the complementary
function to $\widetilde A$.

For $x\in\Omega$ and $s\geq0$, we denote by $a_+^{-1}(x,s)$ the
right-hand derivative of $\widetilde{A}(x,\cdot)$ at $s$ at the same time
define $a_+^{-1}(x,s)= -a_+^{-1}(x,-s)$ for $x\in\Omega$ and
$s\leq0$. Then for $x\in\Omega$ and $s\geq0$, we have
\begin{equation*}
a_+^{-1}(x,s)=\sup\{t\geq0:\, a(x,t)\leq s\}=\inf\{t>0:\,a(x,t)>s\}.
\end{equation*}

\vspace{0.3cm}

\begin{prop}\label{Aa} (See \cite{Fan1, Musielak}.) Let $A\in N(\Omega)$. Then the
following assertions hold:
\begin{enumerate}
\item $A(x,t)\leq a(x,t)t\leq A(x,2t)$ for $x\in\Omega$ and
$t\in\mathbb{R}$;
\item $A$ and $\widetilde A$ satisfy the Young inequality
\begin{equation*}
st\leq A(x,t)+\widetilde A(x,s) \text{ for }x\in\Omega \text{ and }
s,t\in\mathbb{R}
\end{equation*}
and the equality holds if $s=a(x,t)$ or $t=a_+^{-1}(x,s)$.
\end{enumerate}
\end{prop}

\vspace{0.3cm}

Let $A,B\in N(\Omega)$. We say that $A$ is weaker than $B$, denoted
by $A\preccurlyeq B$, if there exist positive constants $K_1,K_2$
and $h\in L^1(\Omega)\cap L^\infty(\Omega)$ such that
\begin{equation}\label{preccurlyeq}
A(x,t)\leq K_1B(x,K_2t)+h(x)\text{ for }x\in\Omega\text{ and
}t\in[0,+\infty).
\end{equation}

\vspace{0.3cm}

\begin{prop}(See \cite{Fan1, Musielak}.)
Let $A, B\in N(\Omega)$ and $A\preccurlyeq B$. Then $\widetilde
B\preccurlyeq \widetilde A$, $L^B(\Omega)\hookrightarrow
L^A(\Omega)$ and $L^{\widetilde A}(\Omega)\hookrightarrow
L^{\widetilde B}(\Omega)$.
\end{prop}

\vspace{0.3cm}

\begin{defi}\label{delta_2}
We say that a function $A:\Omega\times[0,+\infty)\rightarrow[0,+\infty)$
satisfies the $\Delta_2(\Omega)$ condition, denoted by $A\in
\Delta_2(\Omega)$, if there exists a positive constant $K>0$ and a
nonnegative function $h\in L^1(\Omega)$ such that
\begin{equation*}
A(x,2t)\leq KA(x,t)+h(x)\text{ for }x\in\Omega\text{ and
}t\in[0,+\infty).
\end{equation*}
\end{defi}

If $A(x,t)=A(t)$ is an $N$-function and $h(x)\equiv0$ in $\Omega$ in
Definition \ref{delta_2}, then $A\in\Delta_2(\Omega)$ if and only if
$A$ satisfies the well-known $\Delta_2$ condition defined in
\cite{Adams, Trudinger_1971}.

\vspace{0.3cm}

\begin{prop}\label{int_A} (See \cite{Fan1}.)
Let $A\in N(\Omega)$ satisfy $\Delta_2(\Omega)$. Then the following
assertions hold,
\begin{enumerate}
\item $L^A(\Omega)=\{u:\,u \text{ is a measurable function, and
}\int_{\Omega}A(x,|u(x)|)\,\mathrm{d}x<+\infty\}$;
\item $\int_\Omega A(x,|u|)\,\mathrm{d}x<1 \text{ (resp. } =1; >1) \Longleftrightarrow \|u\|_A<1 \text{ (resp. } =1;
>1)$, where $u\in L^A(\Omega)$;
\item $\int_\Omega A(x,|u_n|)\,\mathrm{d}x\rightarrow0 \text{ (resp. } 1; +\infty) \Longleftrightarrow \|u_n\|_A\rightarrow0 \text{ (resp. }1;
+\infty)$, where $\{u_n\}\subset L^A(\Omega)$;
\item $u_n\rightarrow u$ in $L^A(\Omega)\Longrightarrow\int_\Omega \big|A(x,|u_n|)\,\mathrm{d}x-
A(x,|u|)\big|\,\mathrm{d}x\rightarrow0$ as $n\rightarrow\infty$;
\item If $\widetilde{A}$ also satisfies $\Delta_2$, then
\begin{equation*}
\bigg|\int_{\Omega}u(x)v(x)\,\mathrm{d}x\bigg|\leq2\|u\|_A\|v\|_{\widetilde
A},\Forall u\in L^A(\Omega),v\in L^{\widetilde A}(\Omega);
\end{equation*}
\item $a(\cdot,|u(\cdot)|)\in L^{\widetilde A}(\Omega)$ for every $u\in
L^A(\Omega)$.
\end{enumerate}
\end{prop}

\vspace{0.3cm}

The following assumptions will be used.
\vspace{0.3cm}
\begin{enumerate}
\item[$(C_1)$] $\inf_{x\in\Omega}A(x,1)=c_1>0$;
\end{enumerate}

\vspace{0.3cm}

\begin{prop}(See \cite{Fan1}.)
If $A\in N(\Omega)$ satisfies $(C_1)$, then
$L^A(\Omega)\hookrightarrow L^1(\Omega)$ and
$W^{1,A}(\Omega)\hookrightarrow W^{1,1}(\Omega)$.
\end{prop}

\vspace{0.3cm}

Let $A\in  N(\Omega)$ be locally integrable. We will denote
\begin{equation*}
\begin{aligned}
W_{0}^{1,A}(\Omega):&=\overline{C_0^{\infty}(\Omega)}^{\|\,\cdot\,\|_{W^{1,A}(\Omega)}}\\
\mathcal{D}_0^{1,A}(\Omega):&=\overline{C_0^{\infty}(\Omega)}^{\|\nabla\,\cdot\,\|_{L^{A}(\Omega)}}.
\end{aligned}
\end{equation*}
In the case that $\|\nabla u\|_{A}$ is an equivalent norm in
$W_0^{1,A}(\Omega)$,
$W_0^{1,A}(\Omega)=\mathcal{D}_0^{1,A}(\Omega)$.

\vspace{0.3cm}

\begin{prop}(See \cite{Fan1}.)
Let $A\in N(\Omega)$ be locally integrable and satisfy $(C_1)$. Then
\begin{enumerate}
\item the spaces
$W^{1,A}(\Omega),W_{0}^{1,A}(\Omega)$ and
$\mathcal{D}_0^{1,A}(\Omega)$ are separable Banach spaces, and
\begin{equation*}
\begin{aligned}
W_{0}^{1,A}(\Omega)&\hookrightarrow W^{1,A}(\Omega)\hookrightarrow
W^{1,1}(\Omega)\\
\mathcal{D}_0^{1,A}(\Omega)&\hookrightarrow
\mathcal{D}_0^{1,1}(\Omega)=W_{0}^{1,1}(\Omega);
\end{aligned}
\end{equation*}
\item the spaces $W^{1,A}(\Omega),W_{0}^{1,A}(\Omega)$ and
$\mathcal{D}_0^{1,A}(\Omega)$ are reflexive provided $L^A(\Omega)$
is reflexive.
\end{enumerate}
\end{prop}

\vspace{0.3cm}

We always assume that all the Musielak-Orlicz spaces appearing in this paper are reflexive. It is well known that every uniformly convex Banach space is reflexive. For the uniform convexity of Musielak-Orlicz-Sobolev spaces see \cite{Musielak, Fan_Guan}.

\vspace{0.3cm}

\begin{prop}\label{Poincare} (See \cite{Fan1}.)
Let $A,B\in  N(\Omega)$ and $A$ be locally integrable. If there is a
compact imbedding $W^{1,A}(\Omega)\hookrightarrow\hookrightarrow
L^{B}(\Omega)$ and $A\preccurlyeq B$, then there holds the following
Poincar\'{e} inequality
\begin{equation*}
\|u\|_A\leq c\|\nabla u\|_A, \Forall u\in W_0^{1,A}(\Omega),
\end{equation*}
which implies that $\|\nabla\cdot\|_A$ is an equivalent norm in
$W_0^{1,A}(\Omega)$ and
$W_0^{1,A}(\Omega)=\mathcal{D}_0^{1,A}(\Omega)$.
\end{prop}

\vspace{0.3cm}

The following assumptions will be used.

\begin{enumerate}
\item[$(P_1)$] $\Omega\subset\mathbb{R}^n(n\geq2)$ is a bounded domain with the cone property,
and $A\in N(\Omega)$;
\item[$(P_2)$]
$A:\overline{\Omega}\times\mathbb{R}\rightarrow[0,+\infty)$ is
continuous and $A(x,t)\in(0,+\infty)$ for $x\in\overline{\Omega}$
and $t\in(0,+\infty)$.
\end{enumerate}

\vspace{0.3cm}

Let $A$ satisfy $(P_1)$ and $(P_2)$. Denote by $A^{-1}(x,\cdot)$ the
inverse function of $A(x,\cdot)$.  We always assume that the
following condition holds.
\begin{enumerate}
\item[$(P_3)$] $A\in N(\Omega)$ satisfies
\begin{equation}\label{0_1}
\int_0^1\frac{A^{-1}(x,t)}{t^{\frac{n+1}{n}}}\d t<+\infty,\Forall
x\in\overline\Omega.
\end{equation}
\end{enumerate}

\vspace{0.3cm}

Under assumptions $(P_1)$, $(P_2)$ and $(P_3)$, for each
$x\in\overline{\Omega}$, the function
$A(x,\cdot):[0,+\infty)\rightarrow[0,+\infty)$ is a strictly
increasing homeomorphism. Define a function $A_*^{-1}:
\overline{\Omega}\times[0,+\infty)\rightarrow[0,+\infty)$ by
\begin{equation}\label{inversA_*}
A_*^{-1}(x,s)=\int_0^s\frac{A^{-1}(x,\tau)}{\tau^{\frac{n+1}{n}}}\,\mathrm{d}\tau\text{
for }x\in\overline{\Omega} \text{ and }s\in[0,+\infty).
\end{equation}
Then under the assumption $(P_3)$, $A_*^{-1}$ is well defined, and
for each $x\in\overline{\Omega}$, $A_*^{-1}(x,\cdot)$ is strictly
increasing, $A_*^{-1}(x,\cdot)\in C^1((0,+\infty))$ and the function
$A_*^{-1}(x,\cdot)$ is concave.

Set
\begin{equation}\label{T}
T(x)=\lim_{s\rightarrow+\infty}A_*^{-1}(x,s), \Forall
x\in\overline\Omega.
\end{equation}
Then $0<T(x)\leq +\infty$. Define an even function $A_*:
\overline{\Omega}\times\mathbb{R}\rightarrow[0,+\infty)$ by
\begin{equation*}
\begin{aligned}
A_*(x,t)=\left\{ \begin{array}{ll}
          s,  & \text{ if } x\in \overline{\Omega}, |t|\in[0,T(x))\text{ and }A_*^{-1}(x,s)=|t|,\\
          +\infty,   & \text{ for } x\in \overline{\Omega} \text{ and } |t|\geq T(x).
                \end{array}\right.
\end{aligned}
\end{equation*}
Then if $A\in N(\Omega)$ and $T(x)=+\infty$ for any
$x\in\overline{\Omega}$, it is well known that $A_*\in N(\Omega)$
(see \cite{Adams}). $A_*$ is called the Sobolev conjugate function
of $A$ (see \cite{Adams} for the case of Orlicz functions).

Let $X$ be a metric space and $f:X\rightarrow(-\infty,+\infty]$ be
an extended real-valued function. For $x\in X$ with $f(x)\in
\mathbb{R}$, the continuity of $f$ at $x$ is well defined. For $x\in
X$ with $f(x)=+\infty$, we say that $f$ is continuous at $x$ if
given any $M>0$, there exists a neighborhood $U$ of $x$ such that
$f(y)>M$ for all $y\in U$. We say that
$f:X\rightarrow(-\infty,+\infty]$ is continuous on $X$ if $f$ is
continuous at every $x\in X$. Define Dom$(f)=\{x\in X :
f(x)\in\mathbb{R}\}$ and denote by $C^{1-0}(X)$ the set of all
locally Lipschitz continuous real-valued functions defined on $X$.

\vspace{0.3cm}

\begin{rem}\label{rem}
Suppose that $A\in N(\Omega)$ satisfy $(P_2)$. Then for each $t_0\geq0, \widetilde{A}(x,t_0),\\A_*(x,t_0)$ are bounded.
\end{rem}

\vspace{0.3cm}

The following assumptions will also be used.
\begin{enumerate}
\item[$(P_4)$] $T:\overline{\Omega}\rightarrow[0,+\infty]$ is
continuous on $\overline{\Omega}$ and $T\in C^{1-0}(\text{Dom}(T))$;
\item[$(P_5)$]
$A_*\in C^{1-0}(\text{Dom}(A_*))$ and there exist three positive constants
$\delta_0$, $C_0$ and $t_0$ with $\delta_0<\frac{1}{n}$,
$0<t_0<\min_{x\in\overline{\Omega}}T(x)$ such that
\begin{equation*}
|\nabla_x A_*(x,t)|\leq
C_0(A_*(x,t))^{1+\delta_0},
\end{equation*}
for $x\in\Omega$ and $|t|\in[t_0,T(x))$ provided $\nabla_x A_*(x,t)$
exists.
\end{enumerate}

Let $A,B\in N(\Omega)$. We say that $A\ll B$ if, for any $k
> 0$,
\begin{equation*}
\lim_{t\rightarrow+\infty}\frac{A(x,kt)}{B(x,t)}=0\text{ uniformly
for }x\in\Omega.
\end{equation*}

\vspace{0.3cm}

\begin{rem}\label{Symbol}
Suppose that $A,B\in N(\Omega)$. Then $A\ll B\Rightarrow
A\preccurlyeq B$.
\end{rem}

\vspace{0.3cm}

Next we present two embedding theorems for Musielak-Sobolev spaces 
developed by Fan in \cite{Fan2}.

\vspace{0.3cm}

\begin{thm}\label{imbedding}(See \cite{Fan2}, \cite{Liu_Zhao_15}.)
Let $(P_1)-(P_5)$ hold. Then
\begin{enumerate}
\item[(i)] There is a continuous imbedding $W^{1,A}(\Omega)\hookrightarrow
L^{A_*}(\Omega)$;
\item[(ii)] Suppose that $B\in N(\Omega)$,
$B:\overline{\Omega}\times[0,+\infty)\rightarrow[0,+\infty)$ is
continuous, and $B(x,t)\in(0,+\infty)$ for $x\in\Omega$ and
$t\in(0,+\infty)$. If $B\ll A_*$, then there is a compact imbedding
$W^{1,A}(\Omega)\hookrightarrow\hookrightarrow L^B(\Omega)$.
\end{enumerate}
\end{thm}

\vspace{0.3cm}

By Theorem \ref{imbedding}, Remark \ref{Symbol} and Proposition
\ref{Poincare}, we have the following:

\vspace{0.3cm}

\begin{thm}\label{embeddings}(See \cite{Fan2}, \cite{Liu_Zhao_15}.)
Let $(P_1)-(P_5)$ hold and furthermore, $A,A_*\in N(\Omega)$. Then
\begin{enumerate}
\item[(i)] $A\ll A_*$, and there is a compact imbedding $W^{1,A}(\Omega)\hookrightarrow\hookrightarrow
L^A(\Omega)$;
\item[(ii)] there holds the poincar\'{e}-type inequality
\begin{equation*}
\|u\|_A\leq C\|\nabla u\|_A\text{ for }u\in W_0^{1,A}(\Omega),
\end{equation*}
i.e. $\|\nabla u\|_A$ is an equivalent norm on $W_0^{1,A}(\Omega)$.
\end{enumerate}
\end{thm}

\vspace{0.3cm}

In the following of this section, we suppose $A$
satisfies Condition $(\mathscr{A})$, denoted by
$A\in\mathscr{A}$:
\begin{enumerate}
\item[$(\mathscr{A})$]
$A$ satisfies assumptions $(P_1)-(P_3)$, $(P_5)$ in Section \ref{Sec2} and the following
\begin{enumerate}
\item[($\widetilde{P_4}$)] $T(x)$ defined in \eqref{T} satisfies
$T(x)=+\infty$ for all $x\in\overline{\Omega}$.
\end{enumerate}
\end{enumerate}
\vspace{0.3cm}

We recall some useful conclusions in \cite{liu_yao_18} for the reader's convenience.

\vspace{0.3cm}

\begin{lem}\label{compare}(See \cite{liu_yao_18}.)
Suppose that $A\in N(\Omega)$, and there exists a continuous and strictly increasing function
$\mathfrak{A}:[0,+\infty)\rightarrow[0,+\infty)$ such that
\begin{equation}\label{A_A}
A(x,\alpha t)\geq \mathfrak{A}(\alpha)A(x,t), \,\forall\,
\alpha\geq0,t\in\mathbb{R},x\in\Omega.
\end{equation}
\begin{enumerate}
\item[(i)] Then there exists a
continuous and strictly increasing function
$\widehat{\mathfrak{A}}:[0,+\infty)\rightarrow[0,+\infty)$, defined
by
\begin{equation}\label{special}
\begin{aligned}
\widehat{\mathfrak{A}}(\beta)=\left\{ \begin{array}{ll}
          \frac{1}{\mathfrak{A}(\frac{1}{\beta})},  & \text{ for } \beta>0, \\
          0,   & \text{ for } \beta=0,
                \end{array}\right.
\end{aligned}
\end{equation}
such that
\begin{equation}\label{A_star_1}
A(x,\beta t)\leq \widehat{\mathfrak{A}}(\beta)A(x,t), \,\forall
\beta>0,t\in\mathbb{R},x\in\Omega,
\end{equation}
and furthermore $\widehat{\widehat{\mathfrak{A}}}=\mathfrak{A}$;
\item[(ii)] If $\mathfrak{A}$ satisfies
\begin{equation}\label{plessN}
n\mathfrak{A}(\alpha)>\alpha\mathfrak{A}'(\alpha)
\end{equation}
for a.e. $\alpha>0$, then $A_*\in N(\Omega)$, and there exists a
continuous, strictly increasing and a.e. differentiable function
$\mathfrak{A}_*:[0,+\infty)\rightarrow[0,+\infty)$, defined by
\begin{equation}
\begin{aligned}\label{exprr2}
\mathfrak{A}_*^{-1}(\sigma)=\left\{ \begin{array}{ll}
          \frac{1}{\sigma^{\frac{1}{n}}\mathfrak{A}^{-1}(\sigma^{-1})},  & \text{ for } \sigma>0, \\
          0,   & \text{ for } \sigma=0,
                \end{array}\right.
\end{aligned}
\end{equation}
such that
\begin{equation}\label{A_star}
A_*(x,\beta t)\leq \mathfrak{A}_*(\beta)A_*(x,t), \,\forall
\beta>0,t\in\mathbb{R},x\in\Omega;
\end{equation}
\item[(iii)] If $\mathfrak{A}$ satisfies
\begin{equation}\label{plessN2}
\alpha\mathfrak{A}'(\alpha)>\mathfrak{A}(\alpha)
\end{equation}
for a.e. $\alpha>0$, then $\widetilde A\in N(\Omega)$, and there exists a
continuous, strictly increasing and a.e. differentiable function
$\widetilde{\mathfrak{A}}:[0,+\infty)\rightarrow[0,+\infty)$, defined by
\begin{equation}\label{exprr3}
\begin{aligned}
\widetilde{\mathfrak{A}}^{-1}(\sigma)=\left\{ \begin{array}{ll}
         \frac {\sigma}{\mathfrak{A}^{-1}(\sigma)},  & \text{ for } \sigma>0, \\
          0,   & \text{ for } \sigma=0,
                \end{array}\right.
\end{aligned}
\end{equation}
such that
\begin{equation}\label{A_star2}
\widetilde A(x,\beta t)\leq \widetilde{\mathfrak{A}}(\beta)\widetilde A(x,t), \,\forall
\beta>0,t\in\mathbb{R},x\in\Omega.
\end{equation}
\end{enumerate}
\end{lem}

\vspace{0.3cm}

\begin{rem}\label{De_R}
It is clear that $\mathfrak{A}$ defined in \eqref{A_A} depends on $\Omega$. To emphasize the situation we denote $\mathfrak{A}_{\Omega}(t)=\mathfrak{A}(t)$ for any $t\geq0$. To abbreviate the symbols, we write $\mathfrak{A}(t)$ for $\mathfrak{A}_\Omega(t)$ if we consider problems in the domain $\Omega$ in this paper. And it is natural to assume that for $\Omega_0\subset\Omega$
\begin{equation*}
\mathfrak{A}_\Omega(t)\leq\mathfrak{A}_{\Omega_0}(t)\leq\widehat{\mathfrak{A}_{\Omega_0}}(t)\leq\widehat{\mathfrak{A}_\Omega}(t),\,\forall\,t\geq0.
\end{equation*}
\end{rem}

\vspace{0.3cm}

To give more precise increasing assumptions on the function $\mathfrak{A}$ in Lemma \ref{compare} we need the following definition.

\vspace{0.3cm}

\begin{defi}\label{DelT}
\begin{enumerate}
\item Define the operators $\,\,\widehat{\,}$, $\,\,{\,}_*$ and $\,\,\widetilde{\,}$ for any function $\mathfrak{C}:[0,+\infty)\rightarrow[0,+\infty)$ by \eqref{special}, \eqref{exprr2} and \eqref{exprr3} provided $\widehat{\mathfrak{C}}$, ${\mathfrak{C}}_*$ and $\widetilde{\mathfrak{C}}$ exist.
\item We say that $\mathfrak{C}:[0,+\infty)\rightarrow[0,+\infty)$
satisfies Condition $\Delta_{\mathbb{R}^+}$, denoted by
$\mathfrak{C}\in\Delta_{\mathbb{R}^+}$, if there exist positive constants
$M_0$, $M_1$ and $M_2$ such that
\begin{equation}\label{multi}
M_1\mathfrak{C}(\alpha)\mathfrak{C}(\beta)\leq\mathfrak{C}(M_0\alpha\beta)\leq M_2\widehat{\mathfrak{C}}(\alpha)\widehat{\mathfrak{C}}(\beta)
,\,\forall\,\alpha,\beta>0,
\end{equation}
provided $\widehat{\mathfrak{C}}$ exists.
\item We say that $\mathfrak{C}:[0,+\infty)\rightarrow[0,+\infty)$
satisfies Condition $\nabla_{\mathbb{R}^+}$, denoted by
$\mathfrak{C}\in\nabla_{\mathbb{R}^+}$, if there exist positive constants
$M_0$, $M_1$ and $M_2$ such that
\begin{equation}\label{multi1}
M_1\widehat{\mathfrak{C}}(\alpha)\widehat{\mathfrak{C}}(\beta)\leq\mathfrak{C}(M_0\alpha\beta)\leq M_2\mathfrak{C}(\alpha)\mathfrak{C}(\beta)
,\,\forall\,\alpha,\beta>0,
\end{equation}
provided $\widehat{\mathfrak{C}}$ exists.
\end{enumerate}
\end{defi}

\vspace{0.3cm}

\begin{lem}\label{bac_delta1}
Suppose $\mathfrak{C},\mathfrak{D}\in\Delta_{\mathbb{R}^+}$ ($\nabla_{\mathbb{R}^+}$ respectively) and $\widehat{\mathfrak{C}},\widehat{\mathfrak{D}}$ exist. Then 
\begin{enumerate} 
\item[(1)] $\widehat{\mathfrak{C}},\widehat{\mathfrak{D}}\in\nabla_{\mathbb{R}^+}$ ($\Delta_{\mathbb{R}^+}$ respectively);
\item[(2)] for any $C>0$, there exit constants $M=M(C,\mathfrak{C})>0$ and $\widehat M=\widehat M(C,\mathfrak{C})>0$ such that
\begin{equation*}
\begin{aligned}
M\mathfrak{C}(\alpha)\mathfrak{C}(\beta)\leq\mathfrak{C}(C\alpha\beta)&\leq \widehat M\widehat{\mathfrak{C}}(\alpha)\widehat{\mathfrak{C}}(\beta)
,\,\forall\,\alpha,\beta>0,\\
(\widehat M\widehat{\mathfrak{C}}(\alpha)\widehat{\mathfrak{C}}(\beta)\leq\mathfrak{C}(C\alpha\beta)\leq
M&\mathfrak{C}(\alpha)\mathfrak{C}(\beta),\,\forall\,\alpha,\beta>0\text{ respectively}).
\end{aligned}
\end{equation*}
And especially there exit constants $M=M(\mathfrak{C})>0$ and $\widehat M=\widehat M(\mathfrak{C})>0$ such that
\begin{equation*}
\begin{aligned}
M\mathfrak{C}(\alpha)\mathfrak{C}(\beta)\leq\mathfrak{C}(\alpha\beta)&\leq \widehat M\widehat{\mathfrak{C}}(\alpha)\widehat{\mathfrak{C}}(\beta)
,\,\forall\,\alpha,\beta>0,\\
(\widehat M\widehat{\mathfrak{C}}(\alpha)\widehat{\mathfrak{C}}(\beta)\leq\mathfrak{C}(\alpha\beta)\leq
M&\mathfrak{C}(\alpha)\mathfrak{C}(\beta),\,\forall\,\alpha,\beta>0\text{ respectively});
\end{aligned}
\end{equation*}
\item[(3)] $\widehat{\mathfrak{C}}\widehat{\mathfrak{D}}:=\widehat{\mathfrak{C}}(\widehat{\mathfrak{D}}(\cdot))=\widehat{\mathfrak{C}\mathfrak{D}}$;
\item[(4)] $\mathfrak{C}\mathfrak{D}\in\Delta_{\mathbb{R}^+}$ ($\nabla_{\mathbb{R}^+}$ respectively).
\end{enumerate}
\end{lem}

\vspace{0.3cm}

\begin{proof}
(1) Since $\mathfrak{C}\in\Delta_{\mathbb{R}^+}$, there exist constants
$M_0>0$, $M_1>0$ and $M_2>0$ such that for any $\alpha,\beta>0$
\begin{equation*}
M_1\mathfrak{C}(\alpha)\mathfrak{C}(\beta)\leq\mathfrak{C}(M_0\alpha\beta)\leq
M_2\widehat{\mathfrak{C}}(\alpha)\widehat{\mathfrak{C}}(\beta),\,\forall\,\alpha,\beta>0.
\end{equation*}
Then by \eqref{special} we get
\begin{equation*}
\begin{aligned}
\frac{1}{M_2}\mathfrak{C}(\alpha)\mathfrak{C}(\beta)&=\frac{1}{M_2\widehat{\mathfrak{C}}(\alpha^{-1})\widehat{\mathfrak{C}}(\beta^{-1})}\leq\frac{1}{\mathfrak{C}(M_0\alpha^{-1}\beta^{-1})}=
\widehat{\mathfrak{C}}(\frac{1}{M_0}\alpha\beta)\\
&=\frac{1}{\mathfrak{C}(M_0\alpha^{-1}\beta^{-1})}\leq \frac{1}{M_1\mathfrak{C}(\alpha^{-1})\mathfrak{C}(\beta^{-1})}=\frac{1}{M_1}\widehat{\mathfrak{C}}(\alpha)\widehat{\mathfrak{C}}(\beta),
\end{aligned}
\end{equation*}
which implies that $\widehat{\mathfrak{C}}\in\nabla_{\mathbb{R}^+}$.

(2) Since $\mathfrak{C}\in\Delta_{\mathbb{R}^+}$, there exist constants
$M_0>0$, $M_1>0$ and $M_2>0$ such that for any $\alpha,\beta>0$
\begin{equation*}
M_1\mathfrak{C}(\alpha)\mathfrak{C}(\beta)\leq\mathfrak{C}(M_0\alpha\beta)\leq
M_2\widehat{\mathfrak{C}}(\alpha)\widehat{\mathfrak{C}}(\beta),\,\forall\,\alpha,\beta>0.
\end{equation*}
Then by step (1) we get
\begin{equation*}
\begin{aligned}
&M^2_1\mathfrak{C}(\frac{C}{M^2_0})\mathfrak{C}(\alpha)\mathfrak{C}(\beta)\\&\leq M_1\mathfrak{C}(M_0\frac{C}{M^2_0}\alpha)\mathfrak{C}(\beta)=M_1\mathfrak{C}(\frac{C}{M_0}\alpha)\mathfrak{C}(\beta)\leq\mathfrak{C}(M_0\frac{C\alpha}{M_0}\beta)\\&=\mathfrak{C}(C\alpha\beta)\\&=\mathfrak{C}(M_0\frac{C\alpha}{M_0}\beta)\leq
M_2\widehat{\mathfrak{C}}(\frac{C\alpha}{M_1})\widehat{\mathfrak{C}}(\beta)=M_2\widehat{\mathfrak{C}}(\frac{1}{M_0}\frac{CM_0}{M_1}\alpha)\mathfrak{C}(\beta)\\&\leq \frac{M_2}{M_1}\widehat{\mathfrak{C}}(\frac{CM_0}{M_1})\widehat{\mathfrak{C}}(\alpha)\widehat{\mathfrak{C}}(\beta).
\end{aligned}
\end{equation*}
Taking $M:=M^2_1\mathfrak{C}(\frac{C}{M^2_0})$ and $\widehat M:=\frac{M_2}{M_1}\widehat{\mathfrak{C}}(\frac{CM_0}{M_1})$ we get the conclusion of (2).

(3) It is easy to check that for any $\sigma>0$
\begin{equation*}
(\widehat{\mathfrak{C}}\widehat{\mathfrak{D}})(\sigma)=\widehat{\mathfrak{C}}(\widehat{\mathfrak{D}}(\sigma))=\frac{1}{\mathfrak{C}(\frac{1}{\widehat{\mathfrak{D}}(\sigma)})}=\frac{1}{\mathfrak{C}(\mathfrak{D}(\sigma^{-1}))}=\widehat{\mathfrak{C}\mathfrak{D}}(\sigma).
\end{equation*}

(4) It is easy to see that (4) holds. 
\end{proof}

\vspace{0.3cm}

It is easy to see that the following lemma holds.

\vspace{0.3cm}

\begin{lem}\label{bac_delta2}
Suppose $\mathfrak{C}\in\Delta_{\mathbb{R}^+}$ ($\nabla_{\mathbb{R}^+}$ respectively) and $\mathfrak{C}^{-1},\widehat{\mathfrak{C}},\mathfrak{C}_*,\widetilde{\mathfrak{C}}$ exist. Then
\begin{enumerate}
\item $\,\,\widehat{\,}$ is commutative with $\,\,^{-1}$, $\,\,{\,}_*$ and $\,\,\widetilde{\,}$;
\item $\mathfrak{C}^{-1},\widehat{\mathfrak{C}},\mathfrak{C}_*,\widetilde{\mathfrak{C}}\in\nabla_{\mathbb{R}^+}$ ($\Delta_{\mathbb{R}^+}$ respectively) and $\widehat{\mathfrak{C}}^{-1},\mathfrak{C}_*^{-1},\widetilde{\mathfrak{C}}^{-1}\in\Delta_{\mathbb{R}^+}$ ($\nabla_{\mathbb{R}^+}$ respectively).
\end{enumerate}
\end{lem}

\vspace{0.3cm}

\begin{proof}
(1) (a) We will prove $\,\,\widehat{\,}\,$ is commutative with $\,\,^{-1}$. In fact, it is easy to see by \eqref{special} that the following equation holds for any $\sigma>0$,
\begin{equation*}
\widehat{\mathfrak{C}}^{-1}(\sigma)=\frac{1}{\mathfrak{C}^{-1}(\sigma^{-1})}=\widehat{\mathfrak{C}^{-1}}(\sigma).
\end{equation*}

(b) We will prove $\,\,\widehat{\,}\,$ is commutative with $\,\,{\,}_*$. In fact, by \eqref{exprr2} and step (a), for any $\sigma>0$, we get
\begin{equation*}
(\widehat{\mathfrak{C}})_*^{-1}(\sigma)=\frac{1}{\sigma^{\frac{1}{n}}\widehat{\mathfrak{C}}^{-1}(\sigma^{-1})}=\frac{1}{\sigma^{\frac{1}{n}}\widehat{\mathfrak{C}^{-1}}(\sigma^{-1})}=\frac{\mathfrak{C}^{-1}(\sigma)}{\sigma^{\frac{1}{n}}}.
\end{equation*}
On the other hand,
\begin{equation*}
\widehat{\mathfrak{C}_*}^{-1}(\sigma)=\widehat{\mathfrak{C}_*^{-1}}(\sigma)=\frac{1}{\mathfrak{C}_*^{-1}(\sigma^{-1})}=\frac{\mathfrak{C}^{-1}(\sigma)}{\sigma^{\frac{1}{n}}}.
\end{equation*}
Then we conclude $(\widehat{\mathfrak{C}})_*^{-1}(\sigma)=\widehat{\mathfrak{C}_*}^{-1}(\sigma)$.

(c) We will prove $\,\,\widehat{\,}\,$ is commutative with $\,\,\widetilde{\,}$. In fact,
by \eqref{exprr3} and step (a), for any $\sigma>0$, we get
\begin{equation*}
(\widehat{\widetilde{\mathfrak{C}}})^{-1}(\sigma)=(\widehat{\widetilde{\mathfrak{C}}^{-1})}(\sigma)=\frac{1}{{\widetilde{\mathfrak{C}}^{-1}(\sigma^{-1})}}=\frac{1}{\frac{\sigma^{-1}}{\mathfrak{C}^{-1}(\sigma^{-1})}}=\sigma\mathfrak{C}^{-1}(\sigma^{-1}).
\end{equation*}
On the other hand
\begin{equation*}
(\widetilde{\widehat{\mathfrak{C}}})^{-1}(\sigma)=\frac{\sigma}{{\widehat{\mathfrak{C}}^{-1}(\sigma)}}=\frac{\sigma}{{\widehat{\mathfrak{C}^{-1}}(\sigma)}}=\sigma\mathfrak{C}^{-1}(\sigma^{-1}).
\end{equation*}
Then we conclude $(\widehat{\widetilde{\mathfrak{C}}})^{-1}(\sigma)=(\widetilde{\widehat{\mathfrak{C}}})^{-1}(\sigma)$.

(2) Taking $\mathfrak{C}_*$ for example, if $\mathfrak{C}\in\Delta_{\mathbb{R}^+}$, we will prove $\mathfrak{C}_*^{-1}\in\Delta_{\mathbb{R}^+}$ and $\mathfrak{C}_*\in\nabla_{\mathbb{R}^+}$.

Since $\mathfrak{C}\in\Delta_{\mathbb{R}^+}$, there exist constants
$M_0>0$, $M_1>0$ and $M_2>0$ such that for any $\alpha,\beta>0$
\begin{equation}\label{aka}
M_1\mathfrak{C}(\alpha)\mathfrak{C}(\beta)\leq\mathfrak{C}(M_0\alpha\beta)\leq
M_2\widehat{\mathfrak{C}}(\alpha)\widehat{\mathfrak{C}}(\beta),\,\forall\,\alpha,\beta>0.
\end{equation}

(d) Setting $a=\mathfrak{C}(\alpha)$ and $b=\mathfrak{C}(\beta)$ in \eqref{aka}, we get 
\begin{equation*}
\mathfrak{C}^{-1}(M_1ab)\leq
M_0\mathfrak{C}^{-1}(a)\mathfrak{C}^{-1}(b),\,\forall\,a,b>0.
\end{equation*}
Then by \eqref{exprr2} we conclude that for any $\alpha,\beta>0$
\begin{equation*}
\begin{aligned}
\mathfrak{C}_*^{-1}(\frac{1}{M_1}\alpha\beta)&=\frac{M_1^{\frac{1}{n}}}{\alpha^{\frac{1}{n}}\beta^\frac{1}{n}\mathfrak{C}^{-1}(M_1\frac{1}{\alpha\beta})}\geq\frac{M_1^{\frac{1}{n}}}{M_0\alpha^{\frac{1}{n}}\beta^\frac{1}{n}\mathfrak{C}^{-1}(\frac{1}{\alpha})\mathfrak{C}^{-1}(\frac{1}{\beta})}\\
&=\frac{M_1^{\frac{1}{n}}}{M_0}\mathfrak{C}_*^{-1}(\alpha)\mathfrak{C}_*^{-1}(\beta).
\end{aligned}
\end{equation*}

(e) On the other hand, setting $a=\widehat{\mathfrak{C}}(\alpha)$ and $b=\widehat{\mathfrak{C}}(\beta)$ in \eqref{aka}, we get 
\begin{equation*}
\mathfrak{C}^{-1}(M_2ab)\geq
M_0\widehat{\mathfrak{C}}^{-1}(a)\widehat{\mathfrak{C}}^{-1}(b),\,\forall\,a,b>0.
\end{equation*}
Then by \eqref{exprr2} and step (1) we conclude that for any $\alpha,\beta>0$
\begin{equation*}
\begin{aligned}
\mathfrak{C}_*^{-1}(\frac{1}{M_2}\alpha\beta)&=\frac{M_2^{\frac{1}{n}}}{\alpha^{\frac{1}{n}}\beta^\frac{1}{n}\mathfrak{C}^{-1}(M_2\frac{1}{\alpha\beta})}\leq\frac{M_2^{\frac{1}{n}}}{M_0\alpha^{\frac{1}{n}}\beta^\frac{1}{n}\widehat{\mathfrak{C}}^{-1}(\frac{1}{\alpha})\widehat{\mathfrak{C}}^{-1}(\frac{1}{\beta})}\\
&=\frac{M_2^{\frac{1}{n}}}{M_0}(\widehat{\mathfrak{C}})_*^{-1}(\alpha)(\widehat{\mathfrak{C}})_*^{-1}(\beta)=\frac{M_2^{\frac{1}{n}}}{M_0}\widehat{\mathfrak{C}_*}^{-1}(\alpha)\widehat{\mathfrak{C}_*}^{-1}(\beta)\\
&=\frac{M_2^{\frac{1}{n}}}{M_0}\widehat{\mathfrak{C}_*^{-1}}(\alpha)\widehat{\mathfrak{C}_*^{-1}}(\beta).
\end{aligned}
\end{equation*}

(d) and (e) imply that $\mathfrak{C}_*^{-1}\in\Delta_{\mathbb{R}^+}$ and for $a=\mathfrak{C}_*^{-1}(\alpha)$ and $b=\mathfrak{C}_*^{-1}(\beta)$
\begin{equation*}
\mathfrak{C}_*(\frac{M_1^{\frac{1}{n}}}{M_0}ab)\leq\frac{1}{M_1}\mathfrak{C}_*(a)\mathfrak{C}_*(b),
\end{equation*}
and for $a=\widehat{\mathfrak{C}_*}^{-1}(\alpha)$ and $b=\widehat{\mathfrak{C}_*}^{-1}(\beta)$
\begin{equation*}
\mathfrak{C}_*(\frac{M_2^{\frac{1}{n}}}{M_0}ab)\geq\frac{1}{M_2}\widehat{\mathfrak{C}_*}(a)\widehat{\mathfrak{C}_*}(b),
\end{equation*}
which implies $\mathfrak{C}_*\in\nabla_{\mathbb{R}^+}$.

$\widehat{\mathfrak{C}},\widetilde{\mathfrak{C}}\in\nabla_{\mathbb{R}^+}$ and $\widehat{\mathfrak{C}}^{-1},\widetilde{\mathfrak{C}}^{-1}\in\Delta_{\mathbb{R}^+}$ can be obtained by the equation \eqref{special} and equation \eqref{exprr3}.
\end{proof}

\vspace{0.3cm}

Denote the $n$-cube centered at $x_0$ with the edge length of $2R$ by $Q_R(x_0):=\{(x_1,
\dots,x_n)\in\mathbb{R}^n:\,|x_i-(x_0)_i|\leq R, i=1,\dots,n\}$ and 
$u_U=\fint_Uu(x)\,\mathrm{d}x:=|U|^{-1}\int_Uu(x)\,\mathrm{d}x$ for any open set $U
\subset\mathbb{R}^n$ and any $u\in L^1(U)$.

\vspace{0.3cm}
The following two Poincar\'{e} type inequalities in Musielak-Sobolev space are presented with a short proof.
\vspace{0.3cm}

\begin{lem}\label{Poinc} (Poincar\'{e} type inequality.)  
For any bounded, connected, open $U\subset\mathbb{R}^n$ with the cone property, there exists a constant $C=C(n,A,U)$, such that
 \begin{equation*}
 \|u-u_U\|_A\leq C\|\nabla u\|_A,\, \forall u\in W^{1,A}(U),
 \end{equation*}
in which $A\in\mathscr{A}$.
\end{lem}

\vspace{0.3cm}

\begin{proof}
Suppose on the contrary that, for each integer $k\in\mathbb{N}$, there exists a function $u_k\in W^{1,A}(U)$ satisfying 
\begin{equation}\label{contrary_1}
\|u_k-(u_k)_U\|_A\geq k\|\nabla u_k\|_A.
\end{equation}
Define
\begin{equation*}
v_k:=\frac{u_k-(u_k)_U}{\|u_k-(u_k)_U\|_A}\text{ for } k\in\mathbb{N}.
\end{equation*}
Then 
\begin{equation}\label{stadd}
(v_k)_U=0,\,\|v_k\|_A=1,
\end{equation}
and \eqref{contrary_1} implies
\begin{equation}\label{contrary_2}
\|\nabla v_k\|_A<\frac{1}{k}\text{ for } k\in\mathbb{N}.
\end{equation}
It is clear that the functions $\{v_k:k\in\mathbb{N}\}$ are bounded in $W^{1,A}(U)$. By Theorem \ref{embeddings} there exists a subsequence $\{v_{k_j}:j\in\mathbb{N}\}$ of $\{v_k:k\in\mathbb{N}\}$ and a $v\in L^A(U)$ such that $v_{k_j}\rightarrow v$ in $L^A(U)$. Then by \eqref{stadd} it follows that
\begin{equation}\label{contrary_3}
v_U=0,\,\|v\|_A=1.
\end{equation}

On the other hand \eqref{contrary_2} implies for each $i=1,\dots,n$ and any $\phi\in C_0^\infty(\Omega)$ that
\begin{equation*}
\int_Uv\phi_{x_i}\,\mathrm{d}x=\lim_{k_j\rightarrow+\infty}\int_Uv_{k_j}\phi_{x_i}\,\mathrm{d}x=-\lim_{k_j\rightarrow+\infty}\int_U(v_{k_j})_{x_i}\phi\,\mathrm{d}x=0.
\end{equation*}
Then $\nabla v=0$ a.e.. Since $v\in W^{1,A}(U)$ and $U$ is connected, $v$ is a.e. a constant. And this contradicts to \eqref{contrary_3}.
\end{proof}

\vspace{0.3cm}

\begin{lem}\label{Poin_cube} (Poincar\'{e} type inequality for a cube.) 
 There exists a constant $C=C(n,A)$, such that
 \begin{equation}\label{Poinc2}
 \|u-u_R\|_A\leq CR\|\nabla u\|_A,\, \forall u\in W^{1,A}(Q_R(x)),
 \end{equation}
 in which $u_R=u_{Q_R(x)}$ and $A\in\mathscr{A}$.
\end{lem}
\vspace{0.3cm}

\begin{proof}
The case $U=Q_1(0)$ follows from Lemma \ref{Poinc}. In general if $u\in W^{1,A}(Q_R(x))$, denote $v(y):=u(x+ry)$ for any $y\in Q_1(0)$, then $v\in W^{1,A}(Q_1(0))$, which implies
 \begin{equation*}
 \|v-v_1\|_{L^A(Q_1(0))}\leq C\|\nabla v\|_{L^A(Q_1(0))}.
 \end{equation*}
 Changing variables, we get \eqref{Poinc2}.
\end{proof}

\vspace{0.3cm}

\section{Assumptions and lemmas}\label{Sec3}

In this section, in order to get the conclusions of this paper, we state the assumptions on $A$. Firstly, the following assumption on $A$ will be used.

\vspace{0.3cm}

\begin{enumerate}
\item[$(P_0^*)$]  $A\in N(\Omega)$, and for a.e. $x\in\Omega$, $A^{-\frac{1}{n}}(x,t)$ is convex and differentiable on the variable $t$ for $\{t>0\}$.
\end{enumerate}

\vspace{0.3cm}

\begin{defi}\label{defA*}
Under the assumption $(P_0^*)$ we define an even function on the second variable $A^*:\Omega\times\mathbb{R}\rightarrow[0,+\infty)$ as follows:
\begin{equation}\label{deqA*}
(A^*)^{-1}(x,s):=s^\frac{n+1}{n}(A^{-1}(x,s))'_s,
\end{equation}
for a.e. $x\in\Omega$ and $s>0$, where $(A^{-1}(x,s))'_s$ is the derivative of $A^{-1}(x,s)$ on the variable $s$; and $A^*(x,0)=0$ for a.e. $x\in\Omega$.
\end{defi}

\vspace{0.3cm}

\begin{rem}
Under the assumption $(P_0^*)$, $A^*$ in Definition \ref{defA*}  is well defined.
\end{rem}

\vspace{0.3cm}

\begin{proof}
We only need to prove that, for a.e. $x\in\Omega$, the right hand side of \eqref{deqA*} is monotone on $\{s>0\}$. Indeed, setting $A(x,t)=s$, we get
\begin{equation*}
s^\frac{n+1}{n}(A^{-1}(x,s))'_s=\frac{(A(x,t))^{1+\frac{1}{n}}}{A'_t(x,t)}=\frac{1}{-n(A^{-\frac{1}{n}}(x,t))'_t}.
\end{equation*}
Then by the assumption $(P_0^*)$, for a.e. $x\in\Omega$, the right hand side of \eqref{deqA*} is monotone on $\{s>0\}$.
\end{proof}

\vspace{0.3cm}

We give the assumptions on $A^*$ induced by $A$ to proceed.

\vspace{0.3cm}

\begin{enumerate}
\item[$(P^*_1)$] $\Omega\subset\mathbb{R}^n(n\geq2)$ is a bounded domain with the cone property,
and $A^*\in N(\Omega)$;
\item[$(P^*_2)$]
$A^*:\overline{\Omega}\times\mathbb{R}\rightarrow[0,+\infty)$ is
continuous and $A^*(x,t)\in(0,+\infty)$ for $x\in\overline{\Omega}$
and $t\in(0,+\infty)$;

\item[$(P^*_5)$] There exist three positive constants
$\delta_0$, $C_0$ and $t_0$ with $\delta_0<\frac{1}{n}$,
$t_0>0$ such that
\begin{equation*}
|\nabla_x A^*(x,t)|\leq
C_0(A^*(x,t))^{1+\delta_0},\quad{j=1,\dots,n},
\end{equation*}
for $x\in\Omega$ and $|t|\in[t_0,+\infty)$ provided $\nabla_x A^*(x,t)$
exists.
\end{enumerate}

\vspace{0.3cm}

\begin{rem}
Under the assumption $(P^*_0)$, the corresponding assumptions for $A^*$ analogous to $(P_3)$ and $(\widetilde{P_4})$ for $A$ in Section \ref{Sec2} are automatically satisfied.
\end{rem}

\vspace{0.3cm}

Under the assumptions $(P^*_0)$, $(P^*_1)$, $(P^*_2)$ and $(P^*_5)$ it is clear that the following conclusions hold.

\vspace{0.3cm}

\begin{lem}\label{A*sobo}
Assume that $(P^*_0)$, $(P^*_1)$, $(P^*_2)$ and $(P^*_5)$ hold. Then
\begin{enumerate}
\item[(1)] $(A^*)_*=A$;
\item[(2)] $A\in N(\Omega)$ satisfies the assumption $(P_1)$ and $(P_2)$;
\item[(3)] There is a continuous embedding $W^{1,A^*}(\Omega)\hookrightarrow
L^{A}(\Omega)$;
\item[(4)] $A^*\ll A$, and there is a compact embedding $W^{1,A^*}(\Omega)\hookrightarrow\hookrightarrow
L^{A^*}(\Omega)$.
\end{enumerate}
\end{lem}

\vspace{0.3cm}

The following assumption of increasing condition on $A^*\in N(\Omega)$ will be used.

\vspace{0.3cm}

\begin{enumerate}
\item[$(P^*)$] For any $\Omega_0\subseteq\Omega$ with cone property, there exists a strictly increasing differentiable function
$\mathfrak{A}_{\Omega_0}^*\in\nabla_{\mathbb{R}^+}$, with the uniform constants $M_0,M_1$ and $M_2$ for any $\Omega_0\subseteq\Omega$ in Definition \ref{DelT}, such that
\begin{equation}
A^*(x,\alpha t)\leq \mathfrak{A}_{\Omega_0}^*(\alpha)A^*(x,t), \,\forall\,
\alpha\geq0,t\in\mathbb{R},x\in\Omega_0,
\end{equation}
and
\begin{equation}
n\mathfrak{A}_{\Omega_0}^*(\alpha)>\alpha(\mathfrak{A}_{\Omega_0}^*)'(\alpha),\,\forall\,
\alpha\geq0.
\end{equation}
\end{enumerate}

\vspace{0.3cm}

We give a generalized Sobolev Poincar\'{e} type inequality in Musielak-Orlicz-Sobolev space with a short proof.

\vspace{0.3cm}

\begin{lem}\label{SPin} (Sobolev Poincar\'{e} type inequality.) 
For any $Q_R(x)\subset\mathbb{R}^n$, there exists a constant $C=C(n,A^*)$, such that
\begin{equation*}
\widehat{\mathfrak{A}_R}^{-1}(\fint_{Q_R(x)}A(y,\frac{|u-u_R|}{R})\,\mathrm{d}y)\leq C\widehat{\mathfrak{A}_R^*}^{-1}(\fint_{Q_R(x)}A^*(y,|\nabla u|)\,\mathrm{d}y)
\end{equation*}
for any $u\in W^{1,A^*}(Q_R(x))$, where $u_R=u_{Q_R(x)}$ and $\mathfrak{A}_R=\mathfrak{A}_{Q_R(x)}$.
\end{lem}
\vspace{0.3cm}
\begin{proof}
By Lemma \ref{Poin_cube}, for $x=0$ and $R=1$, there exists a constant $C=C(n,A)$ such that
\begin{equation*}
 \|v-v_1\|_{A^*}\leq C\|\nabla v\|_{A^*},\, \forall v\in W^{1,A^*}(Q_1(0)).
\end{equation*} 
Then by Lemma \ref{A*sobo} (3), we get that, for any $v\in W^{1,A^*}(Q_1(0))$
\begin{equation*}
\begin{aligned}
\int_{Q_1(0)}&A(y,|v-v_1|)\,\mathrm{d}y=\int_{Q_1(0)}A(y,\frac{|v-v_1|}{\|v-v_1\|_{A}}\|v-v_1\|_{A})\,\mathrm{d}y\\
&\leq \widehat{\mathfrak{A}_1}(\|v-v_1\|_{A})\int_{Q_1(0)}A(y,\frac{|v-v_1|}{\|v-v_1\|_{A}})\,\mathrm{d}y=\widehat{\mathfrak{A}_1}(\|v-v_1\|_{A})\\
&\leq \widehat{\mathfrak{A}_1}(\|v-v_1\|_{A^*}+\|\nabla v\|_{A^*})\\
&\leq \widehat{\mathfrak{A}_1}(C\|\nabla v\|_{A^*})\\ 
&\leq \widehat{\mathfrak{A}_1}\big(C\widehat{\mathfrak{A}_1^*}^{-1}(\int_{Q_1(0)}\widehat{\mathfrak{A}^*}(\|\nabla v\|_{A^*})A^*(y,\frac{|\nabla v|}{\|\nabla v\|_{A^*}})\,\mathrm{d}y)\big)\\
&\leq \widehat{\mathfrak{A}_1}\big(C\widehat{\mathfrak{A}_1^*}^{-1}(\int_{Q_1(0)}A^*(y,|\nabla v|)\,\mathrm{d}y)\big).\\
\end{aligned}
\end{equation*}
Set $u(x+Ry)=v(y)$ for $y\in Q_1(0)$. Then we get that, for any $u\in W^{1,A^*}(Q_R(x))$
\begin{equation*}
\fint_{Q_R(x)}A(y,|u-u_R|)\,\mathrm{d}y\leq \widehat{\mathfrak{A}_R}\big(C\widehat{\mathfrak{A}_R^*}^{-1}(\fint_{Q_R(x)}A^*(y,R|\nabla u|)\,\mathrm{d}y)\big),
\end{equation*}
or equivalently, by setting $u=\frac{\widetilde u}{R}$, we conclude that, for all $\widetilde u\in W^{1,A^*}(Q_R(x))$,
\begin{equation*}
\widehat{\mathfrak{A}_R}^{-1}(\fint_{Q_R(x)}A(y,\frac{|\widetilde u-\widetilde{u}_R|}{R})\,\mathrm{d}y)\leq C\widehat{\mathfrak{A}_R^*}^{-1}(\fint_{Q_R(x)}A^*(y,|\nabla \widetilde u|)\,\mathrm{d}y),
\end{equation*}
where $C=C(n,A^*)$.
\end{proof}

\vspace{0.3cm}

The following lemma plays an important role in this paper, which is a generalization of Lemma 1 in \cite{Gehring73}.

\vspace{0.3cm}

\begin{lem}\label{ke_ineq}(Gehring type inequality.)
Suppose that $\mathfrak{B}:[0,\infty)\rightarrow[0,\infty)$ is non-decreasing and differentiable, $a\in(1,\infty)$, and $h, H:[\mathfrak{B}^{-1}(1),\infty)\rightarrow[0,\infty)$ are non-increasing with 
\begin{equation}
\lim_{t\rightarrow\infty}h(t)=0,\quad\lim_{t\rightarrow\infty}H(t)=0
\end{equation}
such that for $t\in[\mathfrak{B}^{-1}(1),\infty)$ and a constant $C>0$,
\begin{equation}\label{it_rrr}
-\int_t^\infty\mathfrak{B}(s)\,\mathrm{d}h(s)\leq a\mathfrak{B}(t)h(t)+CH(t).
\end{equation}
Then, for $\epsilon\in[0,\frac{1}{a-1})$,
\begin{equation}\label{Impt_a}
\begin{aligned}
-\int_{\mathfrak{B}^{-1}(1)}^\infty\mathfrak{B}^{1+\epsilon}(t)\,\mathrm{d}h(t)&\leq \frac{1}{1-\epsilon(a-1)}\bigg(-\int_{\mathfrak{B}^{-1}(1)}^\infty\mathfrak{B}(t)\,\mathrm{d}h(t)\bigg)\\
&\quad\quad+\frac{C\epsilon}{1-\epsilon(a-1)}\bigg(-\int_{\mathfrak{B}^{-1}(1)}^\infty\mathfrak{B}^{1+\epsilon}(t)\,\mathrm{d}H(t)\bigg).
\end{aligned}
\end{equation}
\end{lem}

\vspace{0.3cm}

\begin{proof}
\textbf{STEP 1.} Suppose that there exists a $j\in(\mathfrak{B}^{-1}(1),\infty)$ such that $h(t)=0$ and $H(t)=0$ for $t\in[j,\infty)$, and for each $\epsilon\in(0,\infty)$ set
\begin{equation*}
I(\epsilon):=-\int_{\mathfrak{B}^{-1}(1)}^\infty\mathfrak{B}^{1+\epsilon}(t)\,\mathrm{d}h(t)=-\int_{\mathfrak{B}^{-1}(1)}^j\mathfrak{B}^{1+\epsilon}(t)\,\mathrm{d}h(t).
\end{equation*}
Then integration by parts yields
\begin{equation}\label{itttt_1}
\begin{aligned}
I(\epsilon)&=-\int_{\mathfrak{B}^{-1}(1)}^j\mathfrak{B}(t)\mathfrak{B}^{\epsilon}(t)\,\mathrm{d}h(t)\\
&=-\int_{\mathfrak{B}^{-1}(1)}^j\mathfrak{B}(t)\,\mathrm{d}h(t)+\int_{\mathfrak{B}^{-1}(1)}^j\big(1-\mathfrak{B}^{\epsilon}(t)\big)\mathfrak{B}(t)\,\mathrm{d}h(t)\\
&=I(0)+\big(1-\mathfrak{B}^{\epsilon}(t)\big)\bigg(-\int_t^j\mathfrak{B}(s)\,\mathrm{d}h(s)\bigg)\bigg|_{t=\mathfrak{B}^{-1}(1)}^{t=j}\\
&\quad\quad+\epsilon\int_{\mathfrak{B}^{-1}(1)}^j\mathfrak{B}^{\epsilon-1}(t)\bigg(-\int_t^j\mathfrak{B}(s)\,\mathrm{d}h(s)\bigg)\,\mathrm{d}\mathfrak{B}(t)\\
&=I(0)+\epsilon J,
\end{aligned}
\end{equation}
where 
\begin{equation*}
J:=\int_{\mathfrak{B}^{-1}(1)}^j\mathfrak{B}^{\epsilon-1}(t)\bigg(-\int_t^j\mathfrak{B}(s)\,\mathrm{d}h(s)\bigg)\,\mathrm{d}\mathfrak{B}(t).
\end{equation*}
With \eqref{it_rrr} and once more integration by parts we obtain
\begin{equation}\label{itttt_2}
\begin{aligned}
J&\leq a\int_{\mathfrak{B}^{-1}(1)}^j\mathfrak{B}^{\epsilon}(t)h(t)\,\mathrm{d}\mathfrak{B}(t)+C\int_{\mathfrak{B}^{-1}(1)}^j\mathfrak{B}^{\epsilon-1}(t)H(t)\,\mathrm{d}\mathfrak{B}(t)\\
&=\frac{1}{1+\epsilon}[-a\mathfrak{B}(\mathfrak{B}^{-1}(1))h(\mathfrak{B}^{-1}(1))-CH(\mathfrak{B}^{-1}(1))]\\
&\quad\quad-\frac{a}{1+\epsilon}\int_{\mathfrak{B}^{-1}(1)}^j\mathfrak{B}^{1+\epsilon}(t)\,\mathrm{d}h(t)+C\int_{\mathfrak{B}^{-1}(1)}^j\mathfrak{B}^{\epsilon-1}(t)H(t)\,\mathrm{d}\mathfrak{B}(t)\\
&\quad\quad\quad\quad+\frac{C}{1+\epsilon}H(\mathfrak{B}^{-1}(1))\\
&\leq\frac{1}{1+\epsilon}\int_{\mathfrak{B}^{-1}(1)}^\infty\mathfrak{B}(t)\,\mathrm{d}h(t)-\frac{a}{1+\epsilon}\int_{\mathfrak{B}^{-1}(1)}^j\mathfrak{B}^{1+\epsilon}(t)\,\mathrm{d}h(t)\\
&\quad\quad+C\int_{\mathfrak{B}^{-1}(1)}^j\mathfrak{B}^{\epsilon}(t)H(t)\,\mathrm{d}\mathfrak{B}(t)+\frac{C}{1+\epsilon}H(\mathfrak{B}^{-1}(1))\\
&=\frac{1}{1+\epsilon}\int_{\mathfrak{B}^{-1}(1)}^\infty\mathfrak{B}(t)\,\mathrm{d}h(t)-\frac{a}{1+\epsilon}\int_{\mathfrak{B}^{-1}(1)}^j\mathfrak{B}^{1+\epsilon}(t)\,\mathrm{d}h(t)\\
&\quad\quad-\frac{C}{1+\epsilon}\int_{\mathfrak{B}^{-1}(1)}^j\mathfrak{B}^{1+\epsilon}(t)\,\mathrm{d}H(t)\\
&=-\frac{1}{1+\epsilon}I(0)+\frac{a}{1+\epsilon}I(\epsilon)-\frac{C}{1+\epsilon}\int_{\mathfrak{B}^{-1}(1)}^j\mathfrak{B}^{1+\epsilon}(t)\,\mathrm{d}H(t).
\end{aligned}
\end{equation}
Then \eqref{itttt_1} and \eqref{itttt_2} imply that
\begin{equation*}
I(\epsilon)\leq\frac{1}{1-\epsilon(a-1)}I(0)-\frac{C\epsilon}{1-\epsilon(a-1)}\int_{\mathfrak{B}^{-1}(1)}^j\mathfrak{B}^{1+\epsilon}(t)\,\mathrm{d}H(t),
\end{equation*}
whenever $\epsilon\in[0,\frac{1}{a-1})$. So \eqref{Impt_a} follows.

\textbf{STEP 2.} In the general case, for each $j\in(\mathfrak{B}^{-1}(1),\infty)$ set
\begin{equation*}
\begin{aligned}
h_j(t)=\left\{ \begin{array}{ll}
          h(t)  & \text{ if } t\in[\mathfrak{B}^{-1}(1),j), \\
          0   & \text{ if } t\in[j,\infty)
                \end{array}\right.
\end{aligned}
\end{equation*}
and
\begin{equation*}
\begin{aligned}
H_j(t)=\left\{ \begin{array}{ll}
          H(t)  & \text{ if } t\in[\mathfrak{B}^{-1}(1),j), \\
          0   & \text{ if } t\in[j,\infty).
                \end{array}\right.
\end{aligned}
\end{equation*}
Then $h_j,H_j:[\mathfrak{B}^{-1}(1),\infty)\rightarrow[0,\infty)$ is non-increasing and 
for $t\in[\mathfrak{B}^{-1}(1),\infty)$,
\begin{equation*}
-\int_t^\infty\mathfrak{B}(s)\,\mathrm{d}h_j(s)\leq a\mathfrak{B}(t)h_j(t)+CH_j(t).
\end{equation*}
Hence by STEP 1, for $\epsilon\in[0,\frac{1}{a-1})$,
\begin{equation*}
\begin{aligned}
-\int_{\mathfrak{B}^{-1}(1)}^j&\mathfrak{B}^{1+\epsilon}(t)\,\mathrm{d}h(t)=-\int_{\mathfrak{B}^{-1}(1)}^j\mathfrak{B}^{1+\epsilon}(t)\,\mathrm{d}h_j(t)\\
&\leq \frac{1}{1-\epsilon(a-1)}\bigg(-\int_{\mathfrak{B}^{-1}(1)}^j\mathfrak{B}(t)\,\mathrm{d}h_j(t)\bigg)\\
&\quad\quad+\frac{C\epsilon}{1-\epsilon(a-1)}\bigg(-\int_{\mathfrak{B}^{-1}(1)}^j\mathfrak{B}^{1+\epsilon}(t)\,\mathrm{d}H_j(t)\bigg)\\
&\leq \frac{1}{1-\epsilon(a-1)}\bigg(-\int_{\mathfrak{B}^{-1}(1)}^\infty\mathfrak{B}(t)\,\mathrm{d}h_j(t)\bigg)\\
&\quad\quad+\frac{C\epsilon}{1-\epsilon(a-1)}\bigg(-\int_{\mathfrak{B}^{-1}(1)}^\infty\mathfrak{B}^{1+\epsilon}(t)\,\mathrm{d}H_j(t)\bigg).
\end{aligned}
\end{equation*}
Then we obtain \eqref{Impt_a} by letting $j\rightarrow\infty$.
\end{proof}

\vspace{0.3cm}

We get the following main theorem of this paper.

\vspace{0.3cm}

\begin{thm}\label{kkeeyyLemma}
Suppose that $\mathfrak{C}\in\Delta_{\mathbb{R}^+}\cap N $ and $0\leq g\in L^{1}(Q_1(x_0))$ with $\int_{Q_1(x_0)}g(x)\,\mathrm{d}x=1$ and for any $x\in Q_1=Q_1(x_0):=\{x\in\mathbb{R}^n:|x^i-x_0^i|\leq 1,i=1,\dots,n\}$, any $R<\dist(x,\partial Q_1)$, the following estimate (reverse type inequality) holds
\begin{equation} \label{AbL_K_ineq}
\fint_{Q_{\frac{R}{2}}(x)}g\,\mathrm{d}x\leq b\mathfrak{C}(\fint_{Q_R(x)}\mathfrak{C}^{-1}(g)\,\mathrm{d}x)+b
\end{equation}
with $b>0$ being a constant.
Then there exists a constant $\epsilon'>0$, depending only on $\mathfrak{C}$ and $b$, such that for $\epsilon\in[0,\epsilon')$, 
\begin{equation}
g\in L^{\mathfrak{M}}(Q_1),
\end{equation} 
where $\mathfrak{M}(t):=t\cdot(\frac{t}{\mathfrak{C}^{-1}(t)})^\epsilon$; and 
\begin{equation}
\fint_{Q_{\frac{1}{2}}}\mathfrak{M}(g)\,\mathrm{d}x\leq c\mathfrak{M}\bigg(\fint_{Q_1}g\,\mathrm{d}x+1\bigg),
\end{equation}
where $c$ is a positive constant depending only on $n$, $\mathfrak{C}$, $b$ and $\epsilon$. 
\end{thm}

\vspace{0.3cm}

\begin{proof} 
Denote $t_1=t_1(\mathfrak{C})\in(0,+\infty)$ be the unique solution to the equation $\mathfrak{C}(t)=t>0$. Set
\begin{equation*}
\begin{aligned}
C_0&=\{x\in\mathbb{R}^n:|x_i|\leq\frac{1}{2},i=1,\dots,n\}\\
C_k&=\{x\in Q_1:2^{-k-1}\leq\dist(x,\partial Q_1)<2^{-k}\},k\in\mathbb{N}\backslash \{0\}
\end{aligned}
\end{equation*}
and
\begin{equation*}
\begin{aligned}
&G(x)=t_1|Q_1|g(x)\quad 
\text{ on }C_k,\\
&E(G,t)=\{x\in Q_{1}:\,G(x)>t\}.
\end{aligned}
\end{equation*}

\textbf{STEP 1.} We will show that 
\begin{equation}\label{final_equatioooo}
\int_{E(G,t)}G\,\mathrm{d}x\leq a\bigg[\frac{t}{\mathfrak{C}^{-1}(t)}\int_{E(G,t)}\mathfrak{C}^{-1}(G)\,\mathrm{d}x+C\int_{E(F,t)}F\,\mathrm{d}x\bigg]
\end{equation}
for $t\in[t_1,\infty)$, where $a$ is a positive constant depending only on $\mathfrak{C},b$; $C$ is a positive constant depending only on $\mathfrak{C}$ and $F(x):=|Q_1|=\text{a constant}$. 

Fix $t\in[t_1,\infty)$, and for $\epsilon_0>0$, define $s=s(t,\epsilon_0)$ by
\begin{equation*}
s:=bc_3\widehat{\mathfrak{C}}\bigg(\frac{5+(2c_0c_4+1)\epsilon_0}{\epsilon_0c_1^{-1}}\bigg)t,
\end{equation*} 
where $c_0=c_0(\mathfrak{C})$, $c_1=c_1(\mathfrak{C})>0$, $c_3=c_3(\mathfrak{C})>0$ and $c_4=c_4(\mathfrak{C}, g)$ are constants determined in later equations \eqref{Cosss_C_1} and \eqref{apS}. And take $\epsilon_0$ small enough such that $s>t$. It is clear that
\begin{equation*}
Q_1=\cup_{k\geq0}C_k.
\end{equation*}Since 
\begin{equation}\label{t_11111}
\fint_{Q_{1}}G\,\mathrm{d}x=\fint_{Q_1}t_1g(x)|Q_1|\,\mathrm{d}x=t_1,
\end{equation}
we can employ the famous Calder\'{o}n-Zygmund decomposition to obtain a sequence of parallel $n$-cubes $Q^k_j\subset C_{k}$, $j\in\mathbb{N}$, such that
\begin{equation}
G\leq s,\quad  \text{ a.e. in }C_{k}\backslash\cup_{j}Q^k_j
\end{equation}
and
\begin{equation}\label{zhid}
s<\fint_{Q^k_j}G\,\mathrm{d}x\leq 2^ns,\quad\forall j\in\mathbb{N}.
\end{equation}
Then $|E(G,s)\backslash\cup_{k,j}Q^k_j|=0$, which implies from the above two inequalities that
\begin{equation}\label{C_Z_2}
\int_{E(G,s)}G\,\mathrm{d}x\leq\Sigma_{k,j}\int_{Q^k_j}G\,\mathrm{d}x\leq2^ns\Sigma_{k,j}|Q^k_j|.
\end{equation}
Denote by $\overline Q^k_j$ the parallel cube to $Q^k_j$ with the same center and double side, and $D_k:=\cup_j\overline Q^k_j$. It is clear that 
\begin{equation}\label{Set_c}
\overline Q^k_j\subset\cup_{i=k-1}^{i=k+1} C_i.
\end{equation}
Since $D_k$ is bounded, by the well-known covering theorem, there exists a countable disjoint cubes $\{\overline{Q}_{j_h}^k\}_{h=1}^\infty$ which is a subsequence of cubes $\{\overline{Q}_{j}^k\}$ such that
\begin{equation}\label{D_k}
|D_k|\leq5^n\Sigma_h|\overline{Q}_{j_h}^k|.
\end{equation}
Denote by $\overline{Q}_R$ as one of the cubes $\{\overline{Q}_{j_h}^k\}_{h=1}^\infty$. Multiplying \eqref{AbL_K_ineq} by $t_1|Q_1|$, by \eqref{zhid}, we obtain
\begin{equation}\label{Cosss_C_1}
\begin{aligned}
s&<\fint_QG\,\mathrm{d}x\\
&\leq bt_1|Q_1|\mathfrak{C}(\fint_{\overline Q_R}\mathfrak{C}^{-1}(g)\,\mathrm{d}x)+bt_1|Q_1|\\
&\leq \frac{bt_1|Q_1|}{\mathfrak{C}\big(c_1^{-1}\widehat{\mathfrak{C}}^{-1}(t_1|Q_1|)\big)}\mathfrak{C}\bigg(\frac{\widehat{\mathfrak{C}}^{-1}(t_1|Q_1|)}{c_1}\bigg)\mathfrak{C}\bigg(c_1\mathfrak{C}^{-1}(\frac{1}{t_1|Q_1|})\fint_{\overline Q_R}\mathfrak{C}^{-1}(G)\,\mathrm{d}x\bigg)\\
&\quad\quad+bt_1|Q_1|\\
&\leq \frac{bt_1c_2|Q_1|}{\mathfrak{C}\big(c_1^{-1}\widehat{\mathfrak{C}}^{-1}(t_1|Q_1|)\big)}\mathfrak{C}\bigg(\fint_{\overline Q_R}\mathfrak{C}^{-1}(G)\,\mathrm{d}x\bigg)+bt_1|Q_1|\\
&=bc_3\mathfrak{C}\bigg(\fint_{\overline Q_R}\mathfrak{C}^{-1}(G)\,\mathrm{d}x\bigg)+bt_1|Q_1|,
\end{aligned}
\end{equation}
where $c_1=c_1(\mathfrak{C})>0$ satisfies $\mathfrak{C}^{-1}(\alpha\beta)\leq c_1\mathfrak{C}^{-1}(\alpha)\mathfrak{C}^{-1}(\beta),\,\forall\, \alpha,\beta>0$, $c_2=c_2(\mathfrak{C})>0$ satisfies $\mathfrak{C}(\alpha)\mathfrak{C}(\beta)\leq c_2\mathfrak{C}(\alpha\beta),\,\forall\, \alpha,\beta>0$; and $c_3:= \frac{t_1c_2|Q_1|}{\mathfrak{C}\big(c_1^{-1}\widehat{\mathfrak{C}}^{-1}(t_1|Q_1|)\big)}$ Then by the definition of $s$ we can see 
\begin{equation*}
t\leq \mathfrak{C}(\frac{\epsilon_0c_1^{-1}}{5+(2c_0c_4+1)\epsilon_0})\bigg(\mathfrak{C}(\fint_{\overline Q_R}\mathfrak{C}^{-1}(G)\,\mathrm{d}x)+c_3^{-1}t_1|Q_1|\bigg).
\end{equation*}
Since $\mathfrak{C}\in N$ and $\mathfrak{C}^{-1}\in\nabla_{\mathbb{R}^+}$ is concave, it is easy to see that
\begin{equation}\label{c_3}
\mathfrak{C}^{-1}(t)\leq c_1\frac{\epsilon_0c_1^{-1}}{5+(2c_0c_4+1)\epsilon_0}\bigg(\fint_{\overline Q_R}\mathfrak{C}^{-1}(G)\,\mathrm{d}x+\mathfrak{C}^{-1}(c_3^{-1}t_1|Q_1|)\bigg),
\end{equation}
which implies that
\begin{equation}\label{apS}
\begin{aligned}
\frac{5+(2c_0c_4+1)\epsilon_0}{\epsilon_0}\mathfrak{C}^{-1}(t)|\overline Q_R|
&\leq\int_{\overline Q_R}\mathfrak{C}^{-1}(G)\,\mathrm{d}x+\mathfrak{C}^{-1}(c_3^{-1}t_1|Q_1|)|\overline Q_R|\\
&=\int_{\overline Q_R}\mathfrak{C}^{-1}(G)\,\mathrm{d}x+c_4|\overline Q_R|,
\end{aligned}
\end{equation}
where $c_4:=\mathfrak{C}^{-1}(c_3^{-1}t_1|Q_1|)$.
We will estimate the right-hand side of \eqref{apS}. It is easy to see that
\begin{equation}\label{aps_1}
\int_{\overline Q_R}\mathfrak{C}^{-1}(G)\,\mathrm{d}x\leq\int_{E(G,t)\cap\overline Q_R}\mathfrak{C}^{-1}(G)\,\mathrm{d}x+\mathfrak{C}^{-1}(t)|\overline Q_R|.
\end{equation}

At the same time, \eqref{exprr3} and Young inequality for $\mathfrak{C}\in\Delta_{\mathbb{R}^+}$ yield that
\begin{equation}\label{aps_2}
\begin{aligned}
|\overline Q_R|
&=\mathfrak{C}^{-1}(\frac{t|\overline Q_R|}{t})\widetilde{\mathfrak{C}}^{-1}(|\overline Q_R|)\\
&\leq c_0\mathfrak{C}^{-1}(t)\mathfrak{C}^{-1}\bigg(\frac{|\overline Q_R|}{t}\bigg)\widetilde{\mathfrak{C}}^{-1}(|\overline Q_R|)\\
&\leq c_0\mathfrak{C}^{-1}(t)\bigg[\mathfrak{C}^{-1}\bigg(\frac{1}{t}\bigg|E(1,t)\cap \overline Q_R\bigg|+|\overline Q_R|\bigg)\widetilde{\mathfrak{C}}^{-1}(|\overline Q_R|)\bigg]\\
&\leq c_0\frac{\mathfrak{C}^{-1}(t)}{t}\bigg|E(1,t)\cap \overline Q_R\bigg|+2c_0\mathfrak{C}^{-1}(t)|\overline Q_R|.
\end{aligned}
\end{equation}
Therefore, \eqref{apS}, \eqref{aps_1} and \eqref{aps_2} imply that
\begin{equation}\label{wD_k}
\frac{5}{\epsilon_0}|\overline Q_R|\leq \frac{1}{\mathfrak{C}^{-1}(t)}\int_{E(G,t)\cap\overline Q_R}\mathfrak{C}^{-1}(G)\,\mathrm{d}x+c_0c_4\frac{1}{t}\bigg|E(1,t)\cap \overline Q_R\bigg|.
\end{equation}
Combining \eqref{Set_c}, \eqref{D_k} and \eqref{wD_k} we get
\begin{equation*}
\begin{aligned}
|D_k|&\leq\epsilon_05^{n-1}\Sigma_{i=k-1}^{i=k+1}\bigg[  \frac{1}{\mathfrak{C}^{-1}(t)}\int_{E(G,t)\cap C_i}\mathfrak{C}^{-1}(G)\,\mathrm{d}x+c_0c_4\frac{1}{t}\bigg|E(|Q_1|,t)\cap C_i\bigg|\bigg].
\end{aligned}
\end{equation*}
Taking sum with $k$ we can see
\begin{equation*}
\begin{aligned}
\Sigma_k|D_k|\leq\epsilon_0 5^{n}\bigg[  \frac{1}{\mathfrak{C}^{-1}(t)}\int_{E(G,t)}\mathfrak{C}^{-1}(G)\,\mathrm{d}x+c_0c_4\frac{1}{t}\bigg|E(|Q_1|,t)\cap Q_1\bigg|\bigg],
\end{aligned}
\end{equation*}
with \eqref{C_Z_2}, which implies
\begin{equation}\label{FFFana} 
\begin{aligned}
\int_{E(G,s)}G\,\mathrm{d}x\leq a_1b\bigg[\frac{t}{\mathfrak{C}^{-1}(t)}\int_{E(G,t)}\mathfrak{C}^{-1}(G)\,\mathrm{d}x+c_0c_4\bigg|E(|Q_1|,t)\bigg|\bigg]
\end{aligned}
\end{equation}
with 
\begin{equation*}
a_1:=10^n\epsilon_0c_3\widehat{\mathfrak{C}}\bigg(\frac{5+(2c_0c_4+1)\epsilon_0}{\epsilon_0}\bigg).
\end{equation*}
On the other hand, by the definition of $s$ we get
\begin{equation*}
\begin{aligned}
\int_{E(G,t)\backslash E(G,s)}&G\,\mathrm{d}x\leq\int_{E(G,t)}\frac{G}{\mathfrak{C}^{-1}(G)}\mathfrak{C}^{-1}(G)\,\mathrm{d}x\\
&\leq\frac{s}{\mathfrak{C}^{-1}(t)}\int_{E(G,t)}\mathfrak{C}^{-1}(G)\,\mathrm{d}x\\
&\leq bc_3\widehat{\mathfrak{C}}\bigg(\frac{5+(2c_0c_4+1)\epsilon_0}{\epsilon_0}\bigg)\frac{t}{\mathfrak{C}^{-1}(t)}\int_{E(G,t)}\mathfrak{C}^{-1}(G)\,\mathrm{d}x\\
&:= a_2 b\frac{t}{\mathfrak{C}^{-1}(t)}\int_{E(G,t)}\mathfrak{C}^{-1}(G)\,\mathrm{d}x
\end{aligned}
\end{equation*}
with $a_2:=c_3\widehat{\mathfrak{C}}\bigg(\frac{5+(2c_0c_4+1)\epsilon_0}{\epsilon_0}\bigg)$.
The above inequality and \eqref{FFFana} imply
\begin{equation*}
\int_{E(G,t)}G\,\mathrm{d}x\leq a\bigg[\frac{t}{\mathfrak{C}^{-1}(t)}\int_{E(G,t)}\mathfrak{C}^{-1}(G)\,\mathrm{d}x+\frac{c_0c_4}{|Q_1|}\int_{E(F,t)}F\,\mathrm{d}x\bigg]
\end{equation*}
with $a=(a_1+a_2)b$, where $F(x):=|Q_1|=\text{a constant}$. The proof of \eqref{final_equatioooo} is completed.

\textbf{STEP 2.} We will prove the conclusion of the lemma. Now for $t\in[t_1,\infty)$, set $h(t):=\int_{E(G,t)}\mathfrak{C}^{-1}(G)\,\mathrm{d}x=\int_{E(\mathfrak{C}^{-1}(G),\mathfrak{C}^{-1}(t))}\mathfrak{C}^{-1}(G)\,\mathrm{d}x$ and $H(t):=\int_{E(F,t)}F\,\mathrm{d}x=\int_{E(\mathfrak{C}^{-1}(F),\mathfrak{C}^{-1}(t))}\mathfrak{C}(\mathfrak{C}^{-1}(F))\,\mathrm{d}x$. Then $h,H:[t_1,\infty)\rightarrow[0,\infty)$ is non-increasing and
\begin{equation*}
\lim_{t\rightarrow\infty}h(t)=0,\quad\lim_{t\rightarrow\infty}H(t)=0,
\end{equation*} 
and it is easy to verify that
\begin{equation}\label{transf}
\int_{E(P,t)}\mathfrak{D}(P(x))\,\mathrm{d}x=-\int_t^\infty\frac{\mathfrak{D}(s)}{s}\,\mathrm{d}h(s)
\end{equation}
for any $\mathfrak{D}\in N$, any $P\in L^{\mathfrak{D}}(\Omega)$ and any $t\in[t_1(\mathfrak{D}),\infty)$. Then \eqref{final_equatioooo} implies that $h$ satisfies \eqref{it_rrr} for $t'=\mathfrak{C}^{-1}(t)$ by setting $\mathfrak{B}(s)=\frac{\mathfrak{C}(s)}{s}$ for $s\in[t_1,\infty)$ in Lemma \ref{ke_ineq}. Then the conclusion of Lemma \ref{ke_ineq} implies, for $\epsilon\in[0,\frac{1}{a-1})$, that
\begin{equation*}
\begin{aligned}
-\int_{t_1}^\infty\frac{\mathfrak{C}(t')(\frac{\mathfrak{C}(t')}{t'})^\epsilon}{t'}\,\mathrm{d}h(t')&\leq \frac{1}{1-\epsilon(a-1)}\bigg(-\int_{t_1}^\infty\frac{\mathfrak{C}(t')}{t'}\,\mathrm{d}h(t')\bigg)\\
&\quad\quad+\frac{aC\epsilon}{1-\epsilon(a-1)}\bigg(-\int_{t_1}^\infty\frac{\mathfrak{C}(t')(\frac{\mathfrak{C}(t')}{t'})^\epsilon}{t'}\,\mathrm{d}H(t')\bigg).
\end{aligned}
\end{equation*}
Combining with \eqref{transf} we get
\begin{equation*}
\begin{aligned}
\int_{E(G,t_1)}G\cdot(\frac{G}{\mathfrak{C}^{-1}(G)})^\epsilon\,\mathrm{d}x&\leq \frac{1}{1-\epsilon(a-1)}\int_{E(G,t_1)}G\,\mathrm{d}x\\
&\quad\quad+\frac{aC\epsilon}{1-\epsilon(a-1)}\int_{E(F,t_1)}F\cdot(\frac{F}{\mathfrak{C}^{-1}(F)})^\epsilon\,\mathrm{d}x,
\end{aligned}
\end{equation*}
where $C>0$ is a constant depending on $\mathfrak{C}$ and $b$.
Denote $\mathfrak{U}(t):=\big(\frac{t}{\mathfrak{C}^{-1}(t)}\big)^\epsilon$. Since 
\begin{equation*}
G\mathfrak{U}(G)\leq G\text{ in }Q_1\backslash E(G,t_1),
\end{equation*}
we conclude for $\epsilon\in[0,\frac{1}{a-1})$, there exists a constant $C_3>0$ such that
\begin{equation*}
\fint_{Q_1}G\mathfrak{U}(G)\,\mathrm{d}x\leq \frac{C_3}{1-\epsilon(a-1)}\fint_{Q_1}G\,\mathrm{d}x+\frac{C_3\epsilon}{1-\epsilon(a-1)}\fint_{Q_1}F\cdot\mathfrak{U}(F)\,\mathrm{d}x.
\end{equation*}
It is easy to verify $\mathfrak{U}\in\Delta_{\mathbb{R}^+}$. With $\mathfrak{U}(t_1)=1$ and \eqref{t_11111}, we get
\begin{equation*}
\begin{aligned}
\fint_{Q_1}t_1|Q_1|g\cdot\mathfrak{U}(t_1|Q_1|g)\,\mathrm{d}x
&\leq C_4\fint_{Q_1}t_1|Q_1|g\,\mathrm{d}x\cdot\mathfrak{U}\bigg(\fint_{Q_1}t_1|Q_1|g\,\mathrm{d}x\bigg)\\
&\quad\quad+C_4\fint_{Q_1}F\cdot\mathfrak{U}(F)\,\mathrm{d}x,
\end{aligned}
\end{equation*}
which implies that
\begin{equation*}
\fint_{Q_{\frac{1}{2}}}\mathfrak{M}(g)\,\mathrm{d}x\leq c'(n)\fint_{Q_1}\mathfrak{M}(g)\,\mathrm{d}x\leq c\mathfrak{M}\bigg(\fint_{Q_1}g\,\mathrm{d}x+1\bigg),
\end{equation*}
where $c$ is a positive constant depending on $n$, $\mathfrak{C}$, $b$ and $\epsilon$.
\end{proof}

\vspace{0.3cm}

We give another short remark for Theorem \ref{kkeeyyLemma} to end this section.

\vspace{0.3cm}

\begin{rem}
It is easy to verify that $\mathfrak{M}\in N$ and $\mathfrak{M}\gg\mathfrak{C}$.
\end{rem}

\section{The regularity of the minimizers}\label{Sec4}

In the this this section, we suppose $A$
satisfies Condition $\mathscr{A}^*$, denoted by
$A\in\mathscr{A}^*$:
\begin{enumerate}
\item[$(\mathscr{A}^*)$]
$A$ satisfies assumptions $(P^*_0)$, $(P^*_1)$, $(P^*_2)$, $(P^*_5)$ and $(P^*)$ in Section \ref{Sec3}.
\end{enumerate}
\vspace{0.3cm}

Consider the integral functionals as follows
\begin{equation}\label{appf}
E(v)=E(v,\Omega)=\int_\Omega f(x,\nabla
v(x))\,\mathrm{d}x,
\end{equation}
where $v\in W^{1,A}(\Omega)$ and $f(x,z)$ is a Carath\'eodory
function on $\Omega\times\mathbb{R}^n$ satisfying
\begin{equation}\label{appfC}
b_2A\big(x,\sum_{i=1}^n|z_i|\big)-b_1\leq f(x,z)\leq
b_3A\big(x,\sum_{i=1}^n|z_i|\big)+b_1
\end{equation}
with $b_1$, $b_2$ and $b_3$ being non-negative constants, $A\in N(\Omega)\cap\mathscr{A}^*$ satisfying $(C_1)$ (see in Section \ref{Sec2}).

\vspace{0.3cm}

\begin{defi}
A function $u\in W_{loc}^{1,A}(\Omega)$ is said to be a local minimizer of $E$ if
\begin{equation}\label{appmin}
E(u;\text{supp}\varphi)\leq E(u+\varphi;\text{supp}\varphi)\text{ for any }\varphi\in W_0^{1,A}(\Omega)\text{ with } \text{supp}\varphi\subset\subset\Omega.
\end{equation}
\end{defi}

\vspace{0.3cm}

\begin{lem}
Suppose $u\in W_{\text{loc}}^{1,A}$ is a local minimizer of the functional $E$. Then there exist constants $\theta=\theta(b_i,A)\in(0,1)$, $c=c(b_i,A)>0$ and $c'=c'(b_i,A)>0$, such that, for any $x_0\in\Omega$, $Q_R(x_0)\subset\Omega$, $0<t<s<R$ and $a\in(-\infty,+\infty)$, the following inequality holds:
\begin{equation}\label{itera}
\begin{aligned}
\int_{Q_t(x_0)}A(x,|\nabla u|)\,\mathrm{d}x&\leq\theta\int_{Q_s(x_0)}A(x,|\nabla u|)\,\mathrm{d}x\\
&\quad\quad+c\int_{Q_s(x_0)}A(x,\bigg|\frac{u-a}{s-t}\bigg|)\,\mathrm{d}x+c'|Q_s(x_0)|.
\end{aligned}
\end{equation}
\end{lem}
\vspace{0.3cm}
\begin{proof}
For $x_0\in\Omega$, $Q_R(x_0)\subset\Omega$, $0<t<s<R$ and $a\in(-\infty,+\infty)$, choose $\eta\in C_0^\infty(Q_s)$ such that $0\leq\eta\leq1$, $\eta|_{Q_t}\equiv1$ and $|\nabla\eta|\leq\frac{2}{s-t}$. Set $v=u-\eta(u-a)$. Since $u$ is the local minimizer of $E$, we get
\begin{equation}\label{lo}
E(u,Q_s)\leq E(v,Q_s).
\end{equation}
From \eqref{appfC} and \eqref{lo} we conclude
\begin{equation}
-b_1|Q_s|+b_2\int_{Q_s}A(x,|\nabla u|)\,\mathrm{d}x\leq b_1|Q_s|+b_3\int_{Q_s}A(x,|\nabla v|)\,\mathrm{d}x,
\end{equation}
which implies 
\begin{equation}\label{pr_2}
\int_{Q_s}A(x,|\nabla u|)\,\mathrm{d}x\leq b_5\int_{Q_s}A(x,|\nabla v|)\,\mathrm{d}x+b_4|Q_s|,
\end{equation}
where $b_4=\frac{2b_1}{b_2}$ and $b_5=\frac{b_3}{b_2}$. By $\nabla v=(1-\eta)\nabla u-(u-a)\nabla\eta$ we get
\begin{equation}\label{pr_1}
\begin{aligned}
&\int_{Q_s}A(x,|\nabla v|)\,\mathrm{d}x\\&\leq\int_{Q_s}A(x,|(1-\eta)\nabla u|+|(u-a)\nabla\eta|)\,\mathrm{d}x\\
&\leq\int_{Q_s}A(x,2\max\{|(1-\eta)\nabla u|,|(u-a)\nabla\eta|\})\,\mathrm{d}x\\
&\leq\widehat{\mathfrak{A}}(2)\bigg(\int_{Q_s}A(x,|(1-\eta)\nabla u|)\,\mathrm{d}x+\int_{Q_s}A(x,|(u-a)\nabla\eta|)\,\mathrm{d}x\bigg)\\
&\leq\widehat{\mathfrak{A}}(2)\int_{Q_s\backslash Q_t}A(x,|\nabla u|)\,\mathrm{d}x+\big(\widehat{\mathfrak{A}}(2)\big)^2\int_{Q_s}A(x,\bigg|\frac{u-a}{s-t}\bigg|)\,\mathrm{d}x.
\end{aligned}
\end{equation}
Then \eqref{pr_2} and \eqref{pr_1} imply
\begin{equation}\label{pr_3}
\begin{aligned}
&\int_{Q_t}A(x,|\nabla u|)\,\mathrm{d}x\leq\int_{Q_s}A(x,|\nabla u|)\,\mathrm{d}x\\
&\leq c_1\int_{Q_s\backslash Q_t}A(x,|\nabla u|)\,\mathrm{d}x+c_2\int_{Q_s}A(x,\bigg|\frac{u-a}{s-t}\bigg|)\,\mathrm{d}x+b_4|Q_s|,
\end{aligned}
\end{equation}
where $c_1=b_5\widehat{\mathfrak{A}}(2)$ and $c_2=b_5\big(\widehat{\mathfrak{A}}(2)\big)^2$. Adding $c_1\int_{Q_s\backslash Q_t}A(x,|\nabla u|)\,\mathrm{d}x$ on both sides of the \eqref{pr_3}, we conclude
\begin{equation}\label{pr_4}
\begin{aligned}
(1+c_1)\int_{Q_t}A(x,|\nabla u|)\,\mathrm{d}x&\leq c_1\int_{Q_s}A(x,|\nabla u|)\,\mathrm{d}x\\
&\quad\quad+c_2\int_{Q_s}A(x,\bigg|\frac{u-a}{s-t}\bigg|)\,\mathrm{d}x+b_4|Q_s|.
\end{aligned}
\end{equation}
Then we can get the conclusion of the lemma from \eqref{pr_4} if we take $\theta=\frac{c_1}{1+c_1}$, $c=\frac{c_2}{1+c_1}$ and $c'=\frac{b_4}{1+c_1}$.
\end{proof}
\vspace{0.3cm}

\begin{lem}\label{ke_inl}
(Caccioppoli type inequality) Suppose $u\in W_{\text{loc}}^{1,A}(\Omega)$ is a local minimizer of the functional $E$, and there exists a constant $T_0\geq1$ such that $\widehat{\mathfrak{A}}$ induced by $\mathfrak{A}$ satisfies $\widehat{\mathfrak{A}}(T_0)=1$.  Then there exist positive constants $c_3=c_3(b_i,A)$ and $c_4=c_4(b_i,A)$, such that, for any $x_0\in\Omega$, $Q_R(x_0)\subset\Omega$, $0<\rho<R$ and $a\in(-\infty,+\infty)$, the following inequality holds:
\begin{equation}\label{iteee}
\int_{Q_\rho(x_0)}A(x,|\nabla u|)\,\mathrm{d}x\leq c_3\int_{Q_R(x_0)}A(x,\bigg|\frac{u-a}{R-\rho}\bigg|)\,\mathrm{d}x+c_4|Q_R(x_0)|.
\end{equation}
In particular, there exist positive constants $c_5=c_5(b_i, n, A)$ and $c_6=c_6(b_i,n,A)$, such that, for any $x_0\in\Omega$ and $Q_R(x_0)\subset\Omega$, the following inequality holds:
\begin{equation}\label{ittt}
\fint_{Q_{\frac{R}{2}}(x_0)}A(x,|\nabla u|)\,\mathrm{d}x\leq c_5\fint_{Q_R(x_0)}A(x,\bigg|\frac{u-u_R}{R}\bigg|)\,\mathrm{d}x+c_6. 
\end{equation}
\end{lem}

\vspace{0.3cm}

\begin{proof}
Take any $x_0\in\Omega$, $Q_R(x_0)\subset\Omega$, $0<\rho<R$ and $a\in(-\infty,+\infty)$. For $t\in(0,R]$, set $f(t)=\int_{Q_t}A(x,|\nabla u|)\,\mathrm{d}x$. Choose $0<r<1\leq T_0$ such that $\theta\widehat{\mathfrak{A}}(T_0\slash r)<1$. Let
\begin{equation*}
t_0=\rho,\,\, t_{i+1}-t_{i}=\bigg(1-\frac{r}{T_0}\bigg)\bigg(\frac{r}{T_0}\bigg)^i(R-\rho),\,\, i=1,2,\dots,
\end{equation*}
and denote $L:=\int_{Q_R}A(x,\big|\frac{u-a}{R-\rho}\big|)\,\mathrm{d}x$. Iterating \eqref{itera} we get
\begin{equation*}
\begin{aligned}
&f(\rho)=f(t_0)\leq\theta f(t_1)+c\int_{Q_{t_1}}A(x,\bigg|\frac{u-a}{(1-\frac{r}{T_0})(R-\rho)}\bigg|)\,\mathrm{d}x+c'|Q_{t_1}|\\
&\leq\theta f(t_1)+cL\widehat{\mathfrak{A}}(\frac{1}{1-\frac{r}{T_0}})+c'|Q_R|\leq\dots\\
&\leq\theta^k f(t_k)+cL\widehat{\mathfrak{A}}(\frac{1}{1-\frac{r}{T_0}})\Sigma_{i=0}^{k-1}\bigg(\theta\widehat{\mathfrak{A}}(\frac{T_0}{r})\bigg)^i+c'|Q_R|\Sigma_{i=0}^{k-1}\theta^i.
\end{aligned}
\end{equation*}
Sending $k\rightarrow+\infty$ we can get \eqref{iteee}. And \eqref{ittt} is obtained by taking $\rho=R\slash2$ and $a=u_R$ in \eqref{iteee}.
\end{proof}

\vspace{0.3cm}

\begin{lem}\label{endl}
For any $Q_R\subset\Omega$, there exists a positive constant $c_7=c_7(n,A)$ such that, for  
any $u\in W^{1,A}(Q_R)$, the following inequality holds:
\begin{equation}
\fint_{Q_{\frac{R}{2}}}A(x,|\nabla u|)\,\mathrm{d}x\leq c_7\widehat{\mathfrak{A}_R{\mathfrak{A}_R^*}^{-1}}\bigg(\fint_{Q_R}\mathfrak{A}_R^*\mathfrak{A}_R^{-1}(A(x,|\nabla u|))\,\mathrm{d}x\bigg)+c_7.
\end{equation}
\end{lem}

\vspace{0.3cm}

\begin{proof}
Since $A\in N(\Omega)$ satisfies
\begin{equation*}
A(x,\alpha t)\geq \mathfrak{A}_R(\alpha)A(x,t), \,\forall\,
\alpha\geq0,t\in\mathbb{R},x\in Q_R,
\end{equation*}
denoting $\mathfrak{A}_R(\alpha)=\beta$ and $A(x,t)=s$, we can get
\begin{equation*}
A(x,\mathfrak{A}_R^{-1}(\beta)A^{-1}(x,s))\geq \beta s, \,\forall\,
\beta\geq0,s\geq0,x\in\Omega,
\end{equation*}
or equivalently
\begin{equation}\label{kilyd}
A^{-1}(x,\beta s)\leq\mathfrak{A}_R^{-1}(\beta)A^{-1}(x,s) , \,\forall\,
\beta\geq0,s\geq0,x\in\Omega.
\end{equation}
By assumption $(P^*_2)$ and Lemma \ref{A*sobo} (2), there exists a positive constant $\widetilde C$, such that
\begin{equation}\label{CA_ctr}
A^*(x,A^{-1}(x,1))\leq \widetilde C, \,\forall x\in \overline\Omega.
\end{equation}
By $\widehat{\mathfrak{A}}\in\nabla_{\mathbb{R}^+}$, we get
\begin{equation}\label{in_A1111}
\widehat{\mathfrak{A}}(\alpha\beta)\leq C_1\widehat{\mathfrak{A}}(\alpha)\widehat{\mathfrak{A}}(\beta),\,\forall\,\alpha,\beta>0.
\end{equation}
Then by Lemma \ref{ke_inl}, Lemma \ref{SPin}, \eqref{kilyd}, \eqref{CA_ctr},  \eqref{in_A1111} and $\mathfrak{A}^*\in\nabla_{\mathbb{R}_+}$ we get
\begin{equation*}
\begin{aligned}
&\fint_{Q_{\frac{R}{2}}}A(x,|\nabla u|)\,\mathrm{d}x\\
&\leq c_5\fint_{Q_R}A(x,\bigg|\frac{u-u_R}{R}\bigg|)\,\mathrm{d}x+c_6\\
&\leq c_5\widehat{\mathfrak{A}_R}\big(C\widehat{\mathfrak{A}_R^*}^{-1}\bigg(\fint_{Q_R(x)}A^*(x,|\nabla u|)\,\mathrm{d}x\bigg)\big)+c_6\\
&\leq c_5C_1\widehat{\mathfrak{A}_R}(C)\widehat{\mathfrak{A}_R}\widehat{\mathfrak{A}_R^*}^{-1}\bigg(\fint_{Q_R}A^*(x,|\nabla u|)\,\mathrm{d}x\bigg)+c_6\\
&= c_5C_1\widehat{\mathfrak{A}}(C)\widehat{\mathfrak{A}_R}\widehat{\mathfrak{A}_R^*}^{-1}\bigg(\fint_{Q_R}A^*(x,A^{-1}(x,A(x,|\nabla u|)))\,\mathrm{d}x\bigg)+c_6\\
&\leq c_5C_1\widehat{\mathfrak{A}}(C)\widehat{\mathfrak{A}_R}\widehat{\mathfrak{A}_R^*}^{-1}\bigg(\fint_{Q_R}A^*(x,\mathfrak{A}_R^{-1}(A(x,|\nabla u|))A^{-1}(x,1))\,\mathrm{d}x\bigg)+c_6\\
&\leq c_5C_1\widehat{\mathfrak{A}}(C)\widehat{\mathfrak{A}_R}\widehat{\mathfrak{A}_R^*}^{-1}\bigg(\fint_{Q_R}\mathfrak{A}_R^*\mathfrak{A}_R^{-1}(A(x,|\nabla u|))A^*(x,A^{-1}(x,1))\,\mathrm{d}x\bigg)+c_6\\
&\leq c_5C_1\widehat{\mathfrak{A}}(C)\widehat{\mathfrak{A}_R}\widehat{\mathfrak{A}_R^*}^{-1}\bigg(\widetilde C\fint_{Q_R}\mathfrak{A}_R^*\mathfrak{A}_R^{-1}(A(x,|\nabla u|))\,\mathrm{d}x\bigg)+c_6\\
&\leq c_5C_1C_2\widehat{\mathfrak{A}}(C)\widehat{\mathfrak{A}}\widehat{\mathfrak{A}^*}^{-1}(\widetilde C)\widehat{\mathfrak{A}_R}\widehat{\mathfrak{A}_R^*}^{-1}\bigg(\fint_{Q_R}\mathfrak{A}_R^*\mathfrak{A}_R^{-1}(A(x,|\nabla u|))\,\mathrm{d}x\bigg)+c_6\\
&\leq c_7\widehat{\mathfrak{A}_R}\widehat{\mathfrak{A}_R^*}^{-1}\bigg(\fint_{Q_R}\mathfrak{A}_R^*\mathfrak{A}_R^{-1}(A(x,|\nabla u|))\,\mathrm{d}x\bigg)+c_7\\
&= c_7\widehat{\mathfrak{A}_R{\mathfrak{A}_R^*}^{-1}}\bigg(\fint_{Q_R}\mathfrak{A}_R^*\mathfrak{A}_R^{-1}(A(x,|\nabla u|))\,\mathrm{d}x\bigg)+c_7,
\end{aligned}
\end{equation*}
where $c_7=\max\{c_5C_1C_2\widehat{\mathfrak{A}}(C)\widehat{\mathfrak{A}}\widehat{\mathfrak{A}^*}^{-1}(\widetilde C), c_6\}$.
\end{proof}

\vspace{0.3cm}

\begin{lem}\label{imp_co}
Suppose $u\in W^{1,A}(\Omega)$ and there exist two constants $R_0>0$ and $L>0$, not depending on $R$, such that
\begin{equation}\label{inv}
\widehat{\mathfrak{A}_{R}}(t^{-1})\leq L\mathfrak{A}_{R}(t^{-1})
\end{equation}
for any $Q_R\subset Q_1$ and any $t$ with $|Q_R|\leq t\leq |Q_{R_0}|$. Then there exists a constants $c_8=c_8(n, \mathfrak{A}, \Omega, L)$, not depending on $R$, and a constant $R_1$, such that for any $\mathfrak{A}_{R}(|Q_R|)\leq s\leq\mathfrak{A}_{R}(|Q_{R_1}|)$,
\begin{equation}
\widehat{\mathfrak{A}_R{\mathfrak{A}_R^*}^{-1}}(s^{-1})\leq c_8\mathfrak{A}_R{\mathfrak{A}_R^*}^{-1}(s^{-1}).
\end{equation}
\end{lem}

\vspace{0.3cm}

\begin{proof}
Setting $s_{R,t}:=\mathfrak{A}_{R}(t)$, by \eqref{plessN}, there exists a constant $c_0>0$, not depending on $R$, and a constant $R_1\leq R_0$, such that
\begin{equation*}
s_{R,t}^{\frac{1}{n}}=(\mathfrak{A}_{R}(t))^{\frac{1}{n}}\geq  c_0t,\,\forall\,t\text{ with } |Q_R|\leq t\leq| Q_{R_1}|
\end{equation*}
and
\begin{equation*}
s_{R_1,|Q_{R_1}|}^{\frac{1}{n}}=(\mathfrak{A}_{R_1}(|Q_{R_1}|))^{\frac{1}{n}}=c_0|Q_{R_0}|.
\end{equation*}
By \eqref{inv}, we get
\begin{equation}
\widehat{\mathfrak{A}_{R}}(t)\leq L\mathfrak{A}_{R}(t)
\end{equation}
or equivalently,
\begin{equation}
\mathfrak{A}_{R}(t)\mathfrak{A}_{R}(\frac{1}{t})\geq L^{-1}
\end{equation}
for all $t$ with $|Q_R|\leq t\leq |Q_{R_1}|$. Then
\begin{equation*}
\mathfrak{A}_{R}^{-1}(s_{R,t})\leq \widehat{\mathfrak{A}_{R}}^{-1}(Ls_{R,t}).
\end{equation*}
By \eqref{exprr2} it is clear that 
\begin{equation*}
\mathfrak{A}_R^{*-1}(s_{R,t})=\widehat{\mathfrak{A}_R}^{-1}(s_{R,t})s_{R,t}^{\frac{1}{n}}.
\end{equation*}
Then for any $t$ with $|Q_R|\leq t\leq |Q_{R_1}|$, we can see
\begin{equation*}
\begin{aligned}
\widehat{\mathfrak{A}_R\mathfrak{A}_R^{*-1}}(s_{R,t})
&=\frac{1}{\mathfrak{A}_R\bigg(\frac{1}{\mathfrak{A}_R^{-1}(s_{R,t})s_{R,t}^{\frac{1}{n}}}\bigg)}
=\widehat{\mathfrak{A}_R}(\mathfrak{A}_R^{-1}(s_{R,t})s_{R,t}^{\frac{1}{n}})\\
&\leq\widehat{\mathfrak{A}_R}(\widehat{\mathfrak{A}_R}^{-1}(Ls_{R,t})s_{R,t}^{\frac{1}{n}})
\leq C_1Ls_{R,t}\cdot\widehat{\mathfrak{A}_R}(s_{R,t}^{\frac{1}{n}})\\
&=C_1Ls_{R,t}\cdot\widehat{\mathfrak{A}_R}(c_0\frac{s_{R,t}^{\frac{1}{n}}}{c_0})
\leq C_1^2L^2s_{R,t}\cdot\widehat{\mathfrak{A}_R}(c_0)\mathfrak{A}_R(\frac{s_{R,t}^{\frac{1}{n}}}{c_0})\\
&\leq C_1^2L^2\cdot\mathfrak{A}_R\mathfrak{A}_R^{-1}(s_{R,t})
\cdot\frac{\widehat{\mathfrak{A}_R}(c_0)}{\mathfrak{A}_R(c_0)}\mathfrak{A}_R(s_{R,t}^{\frac{1}{n}})\\
&\leq C_2L^2\frac{\widehat{\mathfrak{A}_\Omega}(c_0)}{\mathfrak{A}_\Omega(c_0)}\cdot
\mathfrak{A}_R(\mathfrak{A}_R^{-1}(s_{R,t})s_{R,t}^{\frac{1}{n}})\\
&\leq C_2L^2\frac{\widehat{\mathfrak{A}_\Omega}(c_0)}{\mathfrak{A}_\Omega(c_0)}\cdot
\mathfrak{A}_R(\widehat{\mathfrak{A}_R}^{-1}(s_{R,t})s_{R,t}^{\frac{1}{n}})\\
&:=c_8\mathfrak{A}_R{\mathfrak{A}_R^*}^{-1}(s_{R,t}),
\end{aligned}
\end{equation*}
or equivalently
\begin{equation*}
\begin{aligned}
\widehat{\mathfrak{A}_R\mathfrak{A}_R^{*-1}}(s^{-1})\leq c_8\mathfrak{A}_R\mathfrak{A}_R^{*-1}(s^{-1})
\end{aligned}
\end{equation*}
for any $s$ with $\mathfrak{A}_{R}(|Q_R|)\leq s\leq\mathfrak{A}_{R}(|Q_{R_1}|)$.
\end{proof}

\vspace{0.3cm}

By Jensen's inequality it is easy to get the following lemma.

\vspace{0.3cm}

\begin{lem}\label{est_bdd}
Suppose there exists a constant $R_0>0$, such that for any $Q_R\subset\Omega$ with $R\leq R_0$, $\mathfrak{A}_R\mathfrak{A}_R^{*-1}\in N(Q_R)$.
Then for any $0\leq g\in L^1(\Omega)$, there holds
\begin{equation*} 
\fint_{Q_R}\mathfrak{A}_R^*\mathfrak{A}_R^{-1}(g)\,\mathrm{d}x\leq c_{9}\mathfrak{A}_R^*\mathfrak{A}_R^{-1}\bigg(\fint_{Q_R}g\,\mathrm{d}x\bigg).
\end{equation*}
\end{lem}

\vspace{0.3cm}

\begin{lem}\label{esset}
Suppose $\mathfrak{A}\mathfrak{A}^{*-1}\in N$ and there exist two constants $R_0>0$ and $T_0>0$ such that $\mathfrak{A}^*\mathfrak{A}^{-1}\mathfrak{A}_R\mathfrak{A}_R^{*-1}(t)$ is convex on $\{t\geq T_0\}$ for any $R\leq R_0$.
Then there exists a constant $c_{10}>0$, not depending on $R$, such that for any $0\leq g\in L^1(\Omega)$, there holds
\begin{equation*} 
\mathfrak{A}_R\mathfrak{A}_R^{*-1}\bigg(\fint_{Q_R}\mathfrak{A}_R^*\mathfrak{A}_R^{-1}(g)\,\mathrm{d}x\bigg)\leq c_{10}\mathfrak{A}\mathfrak{A}^{*-1}\bigg(\fint_{Q_R}\mathfrak{A}^*\mathfrak{A}^{-1}(g))\,\mathrm{d}x+1\bigg).
\end{equation*}
\end{lem}

\vspace{0.3cm}

\begin{proof}
It is clear that
\begin{equation*}
\mathfrak{A}_R\mathfrak{A}_R^{*-1}\bigg(\fint_{Q_R}\mathfrak{A}_R^*\mathfrak{A}_R^{-1}(g)\,\mathrm{d}x\bigg)\leq\mathfrak{A}_R\mathfrak{A}_R^{*-1}\bigg(\fint_{Q_R}h(x)\mathrm{d}x\bigg),
\end{equation*}
where 
\begin{equation*}
\begin{aligned}
h(x)=\left\{ \begin{array}{ll}
          \mathfrak{A}_R^*\mathfrak{A}_R^{-1}(g(x)),  & \text{ if } \mathfrak{A}_R^*\mathfrak{A}_R^{-1}(g(x))\geq T_0,\\
          T_0,   & \text{ if } \mathfrak{A}_R^*\mathfrak{A}_R^{-1}(g(x))<T_0.
                \end{array}\right.
\end{aligned}
\end{equation*}
Since $\mathfrak{A}^*\mathfrak{A}^{-1}\mathfrak{A}_R\mathfrak{A}_R^{*-1}(t)$ is convex on $\{t\geq T_0\}$ for any $R\leq R_0$, from Jensen's inequality we get
\begin{equation*}
\begin{aligned}
&\mathfrak{A}_R\mathfrak{A}_R^{*-1}\bigg(\fint_{Q_R}h(x)\mathrm{d}x\bigg)
=\mathfrak{A}\mathfrak{A}^{*-1}\mathfrak{A}^*\mathfrak{A}^{-1}\mathfrak{A}_R\mathfrak{A}_R^{*-1}\bigg(\fint_{Q_R}h(x)\mathrm{d}x\bigg)\\
\leq&\mathfrak{A}\mathfrak{A}^{*-1}\bigg(\fint_{Q_R}\mathfrak{A}^*\mathfrak{A}^{-1}\mathfrak{A}_R\mathfrak{A}_R^{*-1}(h(x))\mathrm{d}x\bigg)\\
\leq&\mathfrak{A}\mathfrak{A}^{*-1}\bigg(\frac{1}{|Q_R|}\int_{\{x\in Q_R:\mathfrak{A}_R^*\mathfrak{A}_R^{-1}(g)\geq T_0\}}\mathfrak{A}^*\mathfrak{A}^{-1}(g)\,\mathrm{d}x\\
&\quad\quad+\frac{|\{x\in Q_R:\mathfrak{A}_R^*\mathfrak{A}_R^{-1}(g)<T_0\}|}{|Q_R|}\mathfrak{A}^*\mathfrak{A}^{-1}\mathfrak{A}_R\mathfrak{A}_R^{*-1}(T_0)\bigg)\\
\leq& \mathfrak{A}\mathfrak{A}^{*-1}\bigg(\fint_{Q_R}\mathfrak{A}^*\mathfrak{A}^{-1}(g)\,\mathrm{d}x+\mathfrak{A}_R\mathfrak{A}_R^{*-1}(T_0)\bigg)\\
\leq& c_{9}\big(\mathfrak{A}\mathfrak{A}^{*-1}\bigg(\fint_{Q_R}\mathfrak{A}^*\mathfrak{A}^{-1}(g)\,\mathrm{d}x\bigg)+\mathfrak{A}_R\mathfrak{A}_R^{*-1}(T_0)\big)\\
\leq& c_{9}\big(\mathfrak{A}\mathfrak{A}^{*-1}\bigg(\fint_{Q_R}\mathfrak{A}^*\mathfrak{A}^{-1}(g)\,\mathrm{d}x\bigg)+\widehat{\mathfrak{A}\mathfrak{A}^{*-1}}(T_0)\big)\\
\leq& c_{10}\big(\mathfrak{A}\mathfrak{A}^{*-1}\bigg(\fint_{Q_R}\mathfrak{A}^*\mathfrak{A}^{-1}(g)\,\mathrm{d}x\bigg)+1\big),
\end{aligned}
\end{equation*}
where $c_9$ and $c_{10}>0$ are constants depending on $\mathfrak{A}$.
\end{proof}

\vspace{0.3cm}

\begin{lem}\label{endlll1}
Under the assumptions of Lemma \ref{imp_co}--\ref{esset}, there exists a constant $c_{14}>0$, not depending on $R$, such that for any $0\leq g\in L^1(\Omega)$ with $\int_{Q_R}g\,\mathrm{d}x\leq1$ and any $Q_R\subset\Omega$,
\begin{equation}
\widehat{\mathfrak{A}_R{\mathfrak{A}_R^*}^{-1}}\bigg(\fint_{Q_R}\mathfrak{A}_R^*\mathfrak{A}_R^{-1}(g)\,\mathrm{d}x\bigg)\leq c_{14}\mathfrak{A}{\mathfrak{A}^*}^{-1}\bigg(\fint_{Q_R}\mathfrak{A}^*\mathfrak{A}^{-1}(g)\,\mathrm{d}x\bigg).
\end{equation}
\end{lem}

\vspace{0.3cm}

\begin{proof}
By Lemma \ref{est_bdd}, \eqref{exprr2} and \eqref{plessN2} we get
\begin{equation*}
\fint_{Q_R}\mathfrak{A}_R^*\mathfrak{A}_R^{-1}(g)\,\mathrm{d}x\leq c_{9}\mathfrak{A}_R^*\mathfrak{A}_R^{-1}\bigg(\frac{1}{|Q_R|}\bigg)\leq c_{12}\frac{1}{|Q_R|}\leq c_{13}\mathfrak{A}_R\bigg(\frac{1}{|Q_R|}\bigg).
\end{equation*}
Then Lemma \ref{imp_co} and Lemma \ref{esset} imply that
\begin{equation*}
\begin{aligned}
\widehat{\mathfrak{A}_R{\mathfrak{A}_R^*}^{-1}}\bigg(\fint_{Q_R}\mathfrak{A}_R^*\mathfrak{A}_R^{-1}(g)\,\mathrm{d}x\bigg)&\leq c_{14}\mathfrak{A}_R\mathfrak{A}_R^{*-1}\bigg(\fint_{Q_R}\mathfrak{A}_R^*\mathfrak{A}_R^{-1}(g)\,\mathrm{d}x+1\bigg)\\
&\leq c_{15}\mathfrak{A}\mathfrak{A}^{*-1}\bigg(\fint_{Q_R}\mathfrak{A}^*\mathfrak{A}^{-1}(g)\,\mathrm{d}x+1\bigg).
\end{aligned}
\end{equation*}
\end{proof}

\vspace{0.3cm}

From Lemma \ref{endl} and Lemma \ref{endlll1} we conclude the following lemma.

\vspace{0.3cm}
\begin{lem}\label{L_4.8}
Under the assumption of Lemma \ref{imp_co}--\ref{esset}, for any $R\leq R_0$, there exists a positive constant $c=c(n,A)$ such that, for  
any $u\in W^{1,A}(Q_R)$, the following inequality holds:
\begin{equation}
\fint_{Q_{\frac{R}{2}}}A(x,|\nabla u|)\,\mathrm{d}x\leq c\mathfrak{A}{\mathfrak{A}^*}^{-1}\bigg(\fint_{Q_R}\mathfrak{A}^*\mathfrak{A}^{-1}(A(x,|\nabla u|))\,\mathrm{d}x\bigg)+c.
\end{equation}
\end{lem}

\vspace{0.3cm}

\begin{lem}\label{LHol}
The following two statement are equivalent:
\begin{enumerate}
\item[(1)] There exists a constant $L>0$, not depending on $R$, such that
\begin{equation}\label{LHCo}
\widehat{\mathfrak{A}_{R}}(|Q_R|^{-1})\leq L\mathfrak{A}_{R}(|Q_R|^{-1})
\end{equation}
for any $Q_R\subset\Omega$.
\item[(2)] There exists a constant $L>0$, not depending on $R$, such that
\begin{equation*}
\widehat{\mathfrak{A}_{R}}(t^{-1})\leq L\mathfrak{A}_{R}(t^{-1})
\end{equation*}
for any $Q_R\subset\Omega$ and any $t\geq|Q_R|$.
\end{enumerate}
\end{lem}

\vspace{0.3cm}

\begin{proof}
If (1) holds, then (2) holds. In fact, from (1), we get
\begin{equation*}
\begin{aligned}
\widehat{\mathfrak{A}_{R}}(t^{-1})\leq \widehat{\mathfrak{A}_{R_t}}(t^{-1})\leq L\mathfrak{A}_{R_t}(t^{-1})\leq \mathfrak{A}_{R}(t^{-1}),
\end{aligned}
\end{equation*}
where $|Q_{R_t}|=t$.
\end{proof}

\vspace{0.3cm}

\begin{rem}\label{L_H_co}
For the variable exponent case $A(x,|t|)=|t|^{p(x)}$, where $1<\inf_{x\in\Omega}p(x)=:p_-\leq p(x)\leq p^+:=\sup_{x\in\Omega}p(x)<n$ and $p\in C(\Omega)$, the log-H\"{o}lder continuity of $p$ can imply that Condition \eqref{LHCo} holds. Check the next section for details.
\end{rem}

Changing variables, by Theorem \ref{kkeeyyLemma}, Lemma \ref{L_4.8} and Lemma \ref{LHol}, we get the following main theorem of this paper.

\vspace{0.3cm}

\begin{thm}\label{main_ltheorem}
Supposee that $(P_0^*)$, $(P_1^*)$, $(P_2^*)$, $(P_5^*)$ and $(P^*)$ hold; that for any $Q_R\subset\Omega$, $\mathfrak{A}_R\mathfrak{A}_R^{*-1}\in N(Q_R)$; that there exists a constant $L>0$, not depending on $R$, such that
\begin{equation*}
\widehat{\mathfrak{A}_{R}}(|Q_R|^{-1})\leq L\mathfrak{A}_{R}(|Q_R|^{-1});
\end{equation*}
that there exist two constants $R_0>0$ and $T_0>0$ such that $\mathfrak{A}^*\mathfrak{A}^{-1}\mathfrak{A}_R\mathfrak{A}_R^{*-1}(t)$ is convex on $\{t\geq T_0\}$ for any $R\leq R_0$. If $u\in W_{loc}^{1,A}(\Omega)$ is a local minimizer of $E$, then there exists a positive constant $\epsilon=\epsilon(n,A)$ such that 
\begin{equation*}
A(x,|\nabla u|)\in L_{\text{loc}}^{\Phi}(\Omega),\quad(\text{equivalently }|\nabla u|\in L_{\text{loc}}^{\Phi(A)}(\Omega),)
\end{equation*} 
where $\Phi(t):=t(\frac{t}{\mathfrak{A}^*{\mathfrak{A}}^{-1}(t)})^\epsilon$ for any $t>0$. And there exists a constant $c=c(n,A,\epsilon,L,|\nabla u|_{L^A(\Omega)})$ such that for every $x\in\Omega$, there exists an $n$-cube $Q_r\subset\Omega$ with the center $x$, such that 
\begin{equation*}
\fint_{Q_{\frac{r}{2}}}\Phi(A(x,|\nabla u|))\,\mathrm{d}x\leq c\Phi\bigg(\fint_{Q_r}A(x,|\nabla u|)\,\mathrm{d}x+1\bigg).
\end{equation*}
\end{thm}

\vspace{0.3cm}

\begin{rem}\label{2}
It is clear that $\Phi\in N$ and $N(\Omega)\ni\Phi(A)\gg A$.
\end{rem}

\section{An Example}\label{Exsec}

\noindent\emph{\textbf{Example.}} Let $p\in
C^{1-0}(\overline\Omega)$ and $\frac{n}{n-1}< q\leq p(x)\leq
p_+:=\text{sup}_{x\in{\overline\Omega}}p(x)<n$ ($q\in \mathbb{R}$)
for $x\in\overline\Omega$. Define
$A:\overline\Omega\times[0,+\infty)\rightarrow[0,+\infty)$ by
\begin{equation*}
A(x,t)=t^{p(x)}.
\end{equation*}
Since $q>\frac{n}{n-1}$,
\begin{equation*}
A^*(x,t)=\big(p(x)t\big)^{\frac{np(x)}{n+p(x)}}\in N(\Omega).
\end{equation*}
It is clear that
\begin{equation*}
A_*(x,t)=\bigg(\frac{n-p(x)}{np(x)}t\bigg)^{\frac{np(x)}{n-p(x)}}\in N(\Omega),
\end{equation*}
and
\begin{equation*}
\begin{aligned}
\mathfrak{A}^*(\alpha):=\left\{ \begin{array}{ll}
          \alpha^{\frac{np^-}{n+p^-}},  & \text{ for } t\geq1, \\
          \alpha^{\frac{np^+}{n+p^+}},  & \text{ for } t<1,
                \end{array}\right.
\end{aligned}
\begin{aligned}
\widehat{\mathfrak{A}^*}(\alpha)=\left\{ \begin{array}{ll}
          \alpha^{\frac{np^+}{n+p^+}},  & \text{ for } t\geq1, \\
          \alpha^{\frac{np^-}{n+p^-}},  & \text{ for } t<1,
                \end{array}\right.
\end{aligned}
\end{equation*}
\begin{equation*}
\begin{aligned}
\mathfrak{A}(\alpha):=\left\{ \begin{array}{ll}
          \alpha^{p^-},  & \text{ for } t\geq1, \\
          \alpha^{p^+},  & \text{ for } t<1,
                \end{array}\right.
\end{aligned}
\begin{aligned}
\widehat{\mathfrak{A}}(\alpha)=\left\{ \begin{array}{ll}
          \alpha^{p^+},  & \text{ for } t\geq1, \\
          \alpha^{p^-},  & \text{ for } t<1
                \end{array}\right.
\end{aligned}
\end{equation*}
\begin{equation*}
\begin{aligned}
\mathfrak{A}_*(\alpha):=\left\{ \begin{array}{ll}
          \alpha^{\frac{np^-}{n-p^-}},  & \text{ for } t\geq1, \\
          \alpha^{\frac{np^+}{n-p^+}},  & \text{ for } t<1,
                \end{array}\right.
\end{aligned}
\begin{aligned}
\widehat{\mathfrak{A}_*}(\alpha)=\left\{ \begin{array}{ll}
          \alpha^{\frac{np^+}{n-p^+}},  & \text{ for } t\geq1, \\
          \alpha^{\frac{np^-}{n-p^-}},  & \text{ for } t<1.
                \end{array}\right.
\end{aligned}
\end{equation*}
The above six control functions satisfy the following inequalities:
\begin{equation*}
\widehat{\mathfrak{A}^*}(\alpha)A(x,t)\leq A^*(x,\alpha t)\leq\mathfrak{A}^*(\alpha)A(x,t),
\end{equation*}
\begin{equation*}
\mathfrak{A}(\alpha)A(x,t)\leq A(x,\alpha t)\leq\widehat{\mathfrak{A}}(\alpha)A(x,t),
\end{equation*}
and
\begin{equation*}
\widehat{\mathfrak{A}_*}(\alpha)A(x,t)\leq A_*(x,\alpha t)\leq\mathfrak{A}_*(\alpha)A(x,t),
\end{equation*}
for $x\in\Omega$ and
$t\in\mathbb{R}^n$. Moreover, we can find that $\mathfrak{A}(\alpha)\in\Delta_{\mathbb{R}_+}$ and $\widehat{\mathfrak{A}^*},\widehat{\mathfrak{A}_*}\in\nabla_{\mathbb{R}_+}$. 

In addition, for $x\in\Omega$,
\begin{equation*}
\nabla_x A^*(x,t)=\bigg[\frac{n^2}{(n+p(x))^2}\log\big(p(x)t\big)+\frac{n}{n+p(x)}\bigg]\big(p(x)t\big)^{\frac{np(x)}{n+p(x)}}\nabla_x p(x).
\end{equation*}
Since for any $\epsilon>0$, $\frac{\log s}{s^\epsilon}\rightarrow0$ as $s\rightarrow+\infty$, we conclude that there
exist constants $\delta_1<\frac{1}{n}$, $c_1$ and $t_1$ such that
\begin{equation*}
\bigg|\frac{\partial A^*(x,t)}{\partial x_j}\bigg|\leq
c_1(A^*)^{1+\delta_1}(x,t),
\end{equation*}
for all $x\in\Omega$ and $t\geq t_1$. Then Condition $(P^*_5)$ is satisfied. 

Denote $p_R:=p|_{Q_R}$. If $p$ is log-H\"{o}lder continuous, i.e. for any $x,y\in Q_R\subset\Omega$ there exists a constant $L>0$ such that
\begin{equation*}
|p_R(x)-p_R(y)|\log|x-y|^{-1}\leq L,
\end{equation*}
then for small $R>0$
\begin{equation*}
R^{-(p_R^+-p_R^-)}\leq|x-y|^{-|p_R(x)-p_R(y)|}\leq \mathrm{e}^L.
\end{equation*}
It is easy to see that
\begin{equation*}
\widehat{\mathfrak{A}_{R}}(|Q_R|^{-1})\leq \mathrm{e}^{nL}\mathfrak{A}_{R}(|Q_R|^{-1}),
\end{equation*}
where $\mathfrak{A}_{R}$ is in the sense of Remark \ref{De_R}.

By the expression of $\mathfrak{A}^*,\mathfrak{A}^*$ and $\mathfrak{A}_*$, the other increasing conditions can be verified. The readers can also find the verification of the increasing conditions in \cite{Fan96} accordingly. This example generalized some parts of the conclusions in \cite{Fan96} and some of the corollaries therein.\\

\section{Acknowledgement}

The authors want to express thanks to Professor Giuseppe Mingione for his encouragement in this research direction and valuable suggestions. The research is supported by the National Natural Science Foundation of China (NSFC 11501268 and NSFC 11471147).
\vspace{0.3cm}
\vspace{0.3cm}
\vspace{0.3cm}



\renewcommand{\baselinestretch}{0.1}
\bibliographystyle{plain}

\bibliography{Ref}

\end{document}